%% file: paper.tex
\documentclass{svjour3}                   
\smartqed  
\usepackage{graphicx}
\usepackage{amsmath}
\usepackage{epstopdf} 
\usepackage{mathptmx}      
\usepackage{booktabs}
\usepackage[T1]{fontenc}
\usepackage[utf8]{inputenc}
\usepackage{graphicx}
\usepackage[ruled,algo2e,linesnumbered,algonl]{algorithm2e}
\usepackage[free-standing-units,group-separator={,}]{siunitx}
\usepackage{caption}
\usepackage{epstopdf}
\usepackage{xcolor}

\graphicspath{{figures/}}

\journalname{BIT}
\begin{document}

\title{New applications for the Boris Spectral Deferred Correction algorithm for plasma simulations\thanks{K.S. was supported by the Engineering and Physical Sciences Research Council (EPSRC) Centre for Doctoral Training in Fluid Dynamics(EP/L01615X/1). SMT would like to acknowledge support from the European Research Council (ERC) under the European Unions Horizon 2020 research and innovation program (grant agreement no. D5S-DLV-786780)}}


\author{Kris Smedt \and
        Daniel Ruprecht \and
        Jitse Niesen \and
        Steven Tobias \and
        Joonas N\"attil\"a
}


\institute{Kris Smedt \at
           Centre for Doctoral Training in Fluid Dynamics, University of Leeds, United Kingdom \\
          \email{kristoffer@smedt.dk}           
           \and
           Daniel Ruprecht \at
           Lehrstuhl Computational Mathematics, Institut für Mathematik, Technische Universität Hamburg, Hamburg, Germany \\
           \email{ruprecht@tuhh.de}
           \and
           Jitse Niesen \at
           School of Mathematics, University of Leeds, United Kingdom \\
           \email{j.niesen@leeds.ac.uk}
            \and
            Steven Tobias \at
            School of Mathematics, University of Leeds, United Kingdom \\
            \email{s.m.tobias@leeds.ac.uk}
            \and
            Joonas N\"attil\"a \at
            Physics Department and Columbia Astrophysics Laboratory, Columbia University, 
            New York, USA \\
            Center for Computational Astrophysics, Flatiron Institute, 
            New York, USA \\
            \email{jan2174@columbia.edu}
}

\date{Received: date / Accepted: date}

\maketitle

\begin{abstract}
The paper investigates two new use cases for the Boris Spectral Deferred Corrections (Boris-SDC) time integrator for plasma simulations.
First, we show that using Boris-SDC as a particle pusher in an electrostatic particle-in-cell (PIC) code can, at least in the linear regime, improve simulation accuracy compared with the standard second order Boris method.
In some instances, the higher order of Boris-SDC even allows a much larger time step, leading to modest computational gains.
Second, we propose a modification of Boris-SDC for the relativistic regime.
Based on an implementation of Boris-SDC in the \textsc{runko} PIC code, we demonstrate for a relativistic Penning trap that Boris-SDC retains its high order of convergence for velocities ranging from $0.5c$ to $>0.99c$.
We also show that for the force-free case where acceleration from electric and magnetic field cancel, Boris-SDC produces less numerical drift than Boris.
\end{abstract}

\input{intro}

\input{methods}

\input{results}

\input{conclusions}

\begin{acknowledgements}
KS thanks NORDITA for the hospitality during his visit during which part of this work was initiated.
\end{acknowledgements}

\bibliographystyle{spmpsci}      
\bibliography{phdSmedt,sdc,additional_refs}

\end{document}

%% file: intro.tex
\section{Introduction}
Movement of charged particles in an electromagnetic field is described by the Lorentz equations
\begin{subequations}
	\begin{align}
	\frac{d \mathbf{x}}{dt}&= \mathbf{v}, \label{eq:lorentz_x} \\
	\frac{d \mathbf{v}}{dt} &= \frac{q}{m} (\mathbf{E} (\mathbf{x}, t) +\mathbf{v} (t) \times\mathbf{B} (\mathbf{x}, t) ) =: \mathbf{f}(\mathbf{x}, \mathbf{v}, t), \label{eq:lorentz_v}
	\end{align}
	\label{eq:lorentz}
\end{subequations}
where $m \mathbf{f}$ is the force on a particle with charge $q$ and mass $m$, $\mathbf{x}$ is the particle position and $\mathbf{v}$ its velocity.
The Lorentz equations have a wide range of applications and are particularly relevant for modeling plasmas.
Understanding plasma dynamics is important since an estimated 99\% of all visible matter in the universe are in a plasma state~\cite{chen1974}.
 
One of the most popular numerical algorithms for solving~\eqref{eq:lorentz} was introduced by Boris in 1970~\cite{boris1970}.
He proposed a Leapfrog method combined with a clever geometric trick to resolve the implicit dependence arising from the $\mathbf{v} \times \mathbf{B}$ term in~\eqref{eq:lorentz_v} based on the observation that the magnetic field only rotates the particle trajectory but does not change the magnitude of its velocity.
Boris' trick can be applied both to leapfrog integration (velocity defined at half time-steps) and Velocity-Verlet integration (velocity and position both defined on integer time-steps). 
In either case, the Boris integrator is second order accurate and conserves phase-space volume, giving it favourable long-term energy behaviour~\cite{qin2013boris}.
There is also a detailed mathematical analysis available, showing that for spatially varying magnetic fields, Boris can still exhibit linear energy drift~\cite{hairer2018energy}.

A number of explicit high order integrators for~\eqref{eq:lorentz} have been developed recently~\cite{hairerLubich2017symmetric,he2016,li2020arbo,qiang2017highOrderRelative,quandt2010high,tao2016explicit}.
Most are derived from the Hamiltonian of~\eqref{eq:lorentz} using splitting methods~\cite{hairer2006geometric}. 
However, very few studies compare them with respect to computational efficiency and none so far investigates their use as a particle pusher in a particle-in-cell code.
Quandt~\cite{quandt2010high} and Li~\cite{li2020arbo} compare computational efficiency of their proposed integrators to classic Boris.
High order integrators were found to show the expected order of convergence and could outperform standard Boris in terms of work-precision for some configurations.

Winkel et al. introduced Boris-SDC in 2015~\cite{winkel2015highOrderBoris}, a combination of the Boris method with the spectral deferred correction (SDC) algorithm by Dutt et al.~\cite{DuttEtAl2000}.
They demonstrate that Boris-SDC delivers high-order accuracy for both a single particle and a particle cloud in a Penning trap and that it leads to less numerical heating than the Boris algorithm.
Tretiak and Ruprecht~\cite{TretiakRuprecht2019}~combine Boris-SDC with a GMRES-based convergence accelerator originally proposed by Huang et al.~\cite{HuangEtAl2006}.
They show that the resulting BGSDC method can deliver improvements in performance over the standard Boris method when simulating fast ions in idealised magnetic fields.
In 2021 they extended these results, showing that BGSDC can improve performance for large ensembles of particles and realistic equilibrium fields of the DIIID and JET Tokamak fusion reactor~\cite{TretiakEtAl2021}.
However, these three  studies consider only non-relativistic cases where particles travel passively through an EM-field guiding them.
Here, we extend their results in two ways. 
First, we investigate numerically the performance of Boris-SDC when used as a particle-pusher in a particle-in-cell code~\cite{hockney1988computer,pukhov2015particle}.
In this case, the particles are no longer passively guided through an electromagnetic field but modify the field.
Second, we introduce a modification to Boris-SDC  for the relativistic case and demonstrate that it retains high order convergence and produces less numerical drift than the Boris method.

While there is some research about potential benefits of using higher order methods in PIC, the focus is mostly on spatial operations like interpolation, deposition or mesh-based approximations of derivatives.
Xiao et al.~\cite{xiao2018geopic} propose a combination of splitting methods applied to the Hamiltonian version of~\eqref{eq:lorentz} and specialised finite elements to produce geometric PIC algorithms capable of high order in both space and time. 
Energy conservation was demonstrated, but no comparison of computational efficiency was made.
Shalaby et al.~\cite{shalaby2017sharp} study the performance of an ESPIC code with higher order algorithms for the field interpolations (up to fifth order) while solving for the field exactly. 
They explicitly highlight the limitation imposed by having only second order time-stepping because of the used Boris/Leapfrog pusher.
Moreover, the energy conservation with charge-conserving PIC algorithms is still being actively studied~\cite{Soklov_2013}.
 
\paragraph{Contributions.}
The paper provides the first investigation of how Boris-SDC performs as a particle pusher in a particle-in-cell code.
It is also one of very few studies that analyses how a high(er) order pusher affects the overall accuracy of PIC.
Furthermore, the paper generalizes the original Boris-SDC algorithm to the relativistic regime.
For three benchmark problems, a two-stream instability, Landau damping and a relativistic Penning trap, we show that the higher order of Boris-SDC leads to substantially better accuracy compared to the standard Boris method and produces less numerical drift in the force-free case, where electric and magnetic field cancel out.
However, in the settings that we tested and without further modifications, the better accuracy is not enough to achieve computational gains from being able to take larger time step: in work-precision studies, the additional work per time step mostly offsets the saved cost from taking fewer but larger steps and Boris-SDC and Boris deliver similar performance.
More details on the results in this paper can be found in the disssertation by Smedt~\cite{Smedt2021}.

%% file: methods.tex
\section{Spectral deferred corrections as pusher for particle-in-cell (PIC)}
Before we describe the Boris spectral deferred correction (Boris-SDC) algorithm, we briefly revisit the key components of the particle-in-cell (PIC) method.
For the sake of simplicity, we restrict our presentation to the one-dimensional case but the generalization to 3D is straightforward, although more complicated in terms of indexing.
A detailed overview is provided for example by Verboncoeur~\cite{verboncoeur2005picReview} and a detailed introduction can be found in the seminal textbook by Birdsall and Landon~\cite{birdsallLangdon1985}.

\subsection{Particle-in-Cell (PIC) \label{sec:2-PIC}}
PIC tracks particles in a Lagrangian approach but has the electric and magnetic field they generate ``live'' on a mesh.
This avoids the $\mathcal{O}(N^2)$ bottleneck that emerges if all particle interactions are computed directly.
The natural consequence of coupling the particles to mesh-based mean fields in PIC is the loss of electrostatic interactions
between particles in close proximity. 
As two particles approach each other, their mutual electric repulsion or attraction goes to zero as opposed to infinity in the real world. 
This emphasises the main assumption of PIC: The global fields arising from the distribution and movement of the full collection of particles are dominant.
Therefore, without additional supplementing algorithms, PIC schemes are only valid for plasmas characterised by collective motion. 

\begin{figure}
    \centering{}
    \includegraphics[width=1\textwidth]{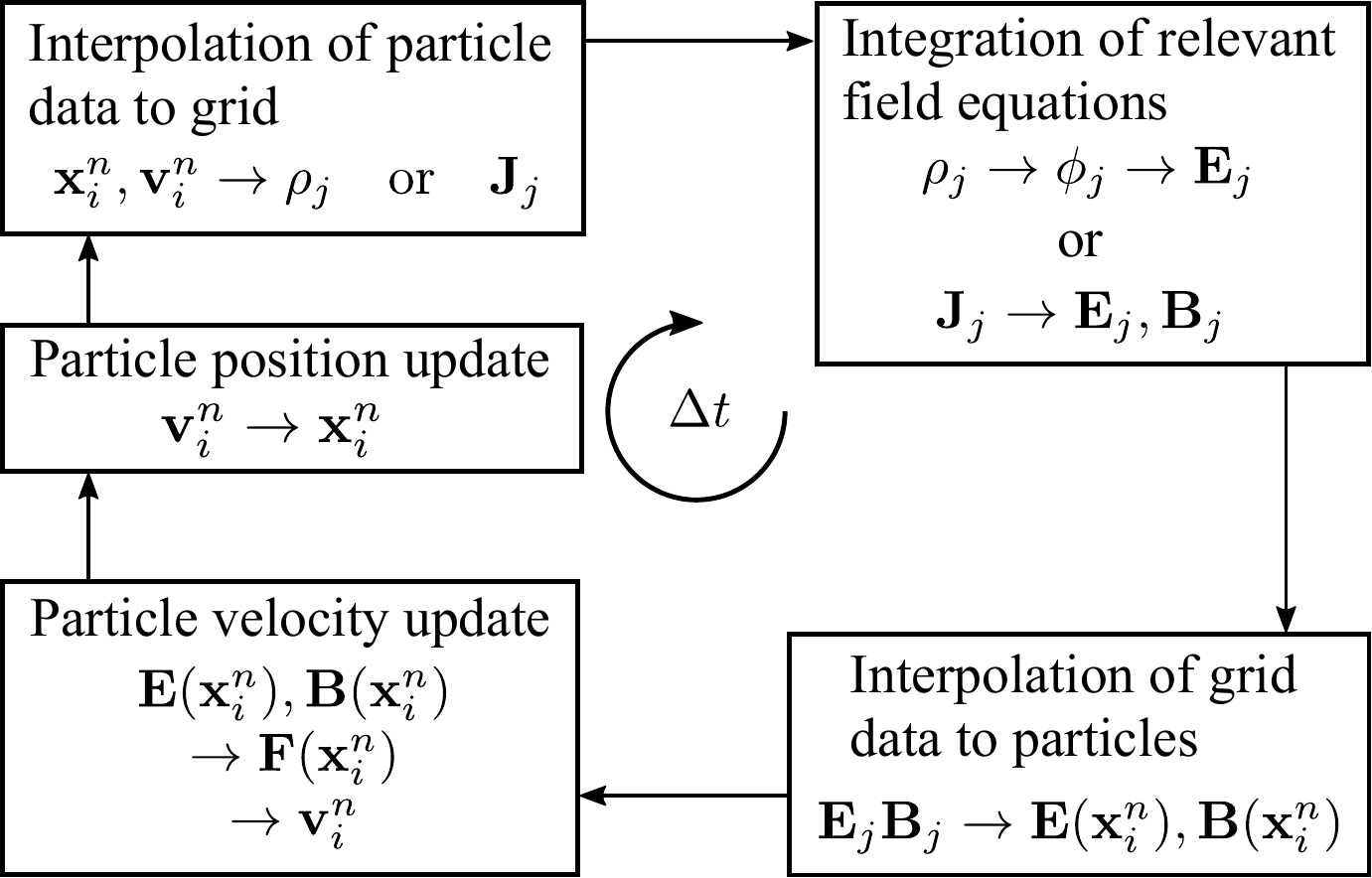} 
    \caption{Steps in a single time step of a general PIC scheme. See e.g. the review by Verboncoeur et al for details~\cite{verboncoeur2005picReview}.}
\label{fig:pic}
\end{figure}

Figure~\ref{fig:pic} sketches the components of one time step in PIC.
This paper focuses on the particle velocity and position update where, using the fields computed in the steps before, particles are moved around by numerically integrating the Lorentz equations.
We investigate how the higher order of accuracy provided by Boris-SDC affects the overall approximation quality of the PIC method.

The exact order of operations in PIC depend on the chosen particle integrator for solving the equation of motion. 
For particle integrators in which the position and velocity are staggered in time, the PIC time-step begins with the velocity update, followed by the position update and field solutions. 
For particle integrators where position and velocity are both defined at the integer time-steps (synchronised), a PIC time-step begins with the position update, followed by the corresponding field solution and finishes with the calculation of the new velocity.
Boris-SDC and Velocity-Verlet, the second order integrator on which it is based, are synchronised particle pushers and so the second type of PIC setup was used throughout this study. 
Whenever the "Boris integrator" term is applied in this study, it refers to Boris' algorithm applied to the velocity-Verlet integrator unless otherwise noted. 
However, note that owing to the popularity of the leapfrog integrator in PIC, most existing schemes are of the staggered type.

\subsubsection{Integration of relevant field equations}
In the electromagnetic case, time derivatives for the electric and magnetic field are present in the Maxwell-Vlasov equations that need to be integrated numerically.
While this step is shown in Fig.~\ref{fig:pic} for the sake of completeness, we only study electrostatic examples in this paper where the electric field is fully reconstructed from the particle charges in every time step.
A detailed survey of different numerical approaches to electromagnetic PIC is provided by Birdsall and Langdon~\cite{birdsallLangdon1985}.
In particular, we do not discuss the issues of divergence correction or cleaning that arises if Gauss' law is not exactly satsified on the discrete level~\cite{MunzEtAl2000} and leave this for future work.

\subsubsection{Interpolation of grid data to particles}\label{subsubsec:interpolation}
In the electrostatic case, the electric field is given by
\begin{equation}
    \mathbf{E} = \nabla \phi
\end{equation}
where $\phi$ is the electrostatic potential.
The potential depends on the charge density $\rho$ via
\begin{equation}
\nabla^{2}\phi (\mathbf{x}, t) =\frac{\rho (\mathbf{x}, t)}{\epsilon}, \label{eq:poisson}
\end{equation}
using the permittivity of the plasma $\epsilon$. Now consider $N_{p}$ particles, where $n = 1, \ldots, N_p$ labels a given particle with position vector $\mathbf{x}_n$ in continuous space and velocity $\mathbf{v}_n$.
Let $x_i$ with $i = 1, \ldots, N_x$ be a set of equidistant mesh points with spacing $x_{i+1} - x_i = \Delta x$.
To calculate the electric field on the grid, we need to determine the corresponding charge densities $\rho_{i}$ that are generated by the particles.
To do this, particle charges $q$ are interpolated to the surrounding grid nodes via some weighting function $W (\mathbf{x}_{n})$. Here, we use linear weighting to ``scatter'' the particle charge to the two nearest grid nodes so that
\begin{equation}
W (x) =\left\{ \begin{array}{cc}
1-\frac{| x - x_i|}{\Delta X} & \, \, \, \, | x - x_{i}|<\Delta x\\
0 & \, \, \, \, |x-x_{i}|>\Delta x.
\end{array}\right. \label{eq:wLinear}
\end{equation}
Other interpolation schemes exist, such as Nearest-Grid-Point or higher order quadratic or cubic weighting splines, but the linear scheme is most commonly used~\cite{verboncoeur2005picReview}. 
Once the charge of the surrounding volume has been assigned to a grid node, the charge density $\rho_{i}$ for the cell is computed as the average over the cell volume. 
Knowledge of the charge densities $\rho_{i}$  allows the determination of  the electrostatic potential with an appropriate solution scheme. For the current study, 1D second order central finite difference was used for the grid.
Discretizing~\eqref{eq:poisson} with second order centered finite differences yields
\begin{equation}
	\frac{\phi_{i-1}-2\phi_{i}+\phi_{i+1}}{\Delta x^2}=\frac{\rho_{i}}{e}.
\end{equation}
The resulting linear system is solved for the $\phi_i$ using the SciPy linear algebra package~\cite{scipy2020}. 
From the grid values, the gradient of the electric potential $\phi_{i}$ and thus the electric field values are computed with a central difference scheme
\begin{equation}
\mathbf{E}_i = \frac{\phi_{i-1}-\phi_{i+1}}{2\Delta x}.
\end{equation}
At boundary nodes, forward
\begin{equation}
\mathbf{E}_1 = \frac{\phi_{1}-\phi_{2}}{\Delta x},
\end{equation}
or backward finite differences 
\begin{equation}
\mathbf{E}_{N_x} = \frac{\phi_{N_x}-\phi_{N_x - 1}}{\Delta x},
\end{equation}
are used instead.

\subsubsection{Particle velocity and position update}
The Newton-Lorentz force gives the acceleration exerted on the particles.
The corresponding differential equation~\eqref{eq:lorentz} is integrated numerically to update velocity and position of the particles from time $t_n$ to time $t_{n} + \Delta t = t_{n+1}$.
A popular algorithm is the St\"ormer-Verlet scheme
\begin{subequations}
\begin{align}
	\mathbf{v}_{n+1/2} &= \mathbf{v}_{n-1/2} + \frac{\Delta t}{2} \mathbf{f}(\mathbf{x}_n, \mathbf{v}_n) \label{eq:stoermer_verlet_v} \\
	\mathbf{x}_{n+1} &= \mathbf{x}_n + \Delta t \mathbf{v}_{n+1/2}
\end{align}	
\end{subequations}
where the calculation of position and velocity are offset by $\Delta t/2$.
Here, $\mathbf{x}_{n} \approx \mathbf{x} (t) $ and $\mathbf{x}_{n+1} \approx \mathbf{x} (t+\Delta t) $, etc. 
Note that some form of interpolation is required to provide $\textbf{v}_n$ in~\eqref{eq:stoermer_verlet_v}.
Typically, the average of $\mathbf{v}_{n+1/2}$ and $\mathbf{v}_{n-1/2}$ is used.
This staggering is advantageous especially on fully electromagnetic PIC loops where the electromagnetic fields can then be evolved with a finite-difference time domain (FDTD) method relying on the so-called Yee lattice \cite{verboncoeur2005picReview}.
However, for Boris-SDC, staggering was found to increase storage requirements without adding much benefit~\cite{TretiakRuprecht2019}. 
We therefore use the second order accurate velocity-Verlet scheme 
\begin{subequations}
	\begin{align}
		\mathbf{x}_{n+1} &= \mathbf{x}_{n}+\mathbf{v}_{n}\Delta t+\frac{1}{2}\mathbf{f} (\mathbf{x}_{n}, \mathbf{v}_{n}) \Delta t^{2}, \label{eq:vv_pos} \\
		\mathbf{v}_{n+1} &= \mathbf{v}_{n}+\frac{\mathbf{f} (\mathbf{x}_{n}, \mathbf{v}_{n}) +\mathbf{f} (\mathbf{x}_{n+1}, \mathbf{v}_{n+1}) }{2}\Delta t, \label{eq:vv_vel}
	\end{align}
	\label{eq:vv}
\end{subequations}
instead.
Both variants are second order accurate and behave very similarly, but they are not equivalent~\cite{Mazur1997}.
Boris-SDC, introduced in detail below, is a high-order generalization of~\eqref{eq:vv}.

\paragraph{Boris' trick.}
\begin{algorithm2e}[t]
        \caption{Boris' trick as a general solver for~\eqref{eq:boris_trick_form}. See Birdsall and Langdon~\cite[Section 4--4]{birdsallLangdon1985} for the geometric derivation.}\label{alg:boris}
        \SetKwComment{Comment}{\# }{}
        \SetCommentSty{textit}
        \SetKwInOut{Input}{input}
    \SetKwInOut{Output}{output}
    \Input{$\mathbf{v}_{n-1}$, $\alpha$, $\beta$, $\mathbf{B}$, $\mathbf{E}$, $\mathbf{c}$}
    \Output{$\mathbf{v}_{n}$ solving  $\mathbf{v}_n = \mathbf{v}_{n-1} + \alpha \mathbf{E} + \beta \frac{\mathbf{v}_{n-1} + \mathbf{v}_n}{2} \times \mathbf{B} + \mathbf{c}$}
    $\mathbf{t} = \frac{\beta}{2} \mathbf{B}$ \\
    $\mathbf{s} = 2 \mathbf{t} / \left( 1 + \mathbf{t} \cdot \mathbf{t} \right)$ \\
    $\mathbf{v}^{-} = \mathbf{v}_{n-1} + \frac{\alpha}{2} \mathbf{E} + \frac{1}{2} \mathbf{c}$ \\
    $\mathbf{v}^{*} = \mathbf{v}^{-} + \mathbf{v}^{-} \times \mathbf{t}$ \\
    $\mathbf{v}^{+} = \mathbf{v}^{-} + \mathbf{v}^{*} \times \mathbf{s}$ \\
    $\mathbf{v}_{n} = \mathbf{v}^{+} + \frac{\alpha}{2} \mathbf{E} + \frac{1}{2} \mathbf{c}$
\end{algorithm2e}
While the position update~\eqref{eq:vv_pos} is explicit, the update for the velocity~\eqref{eq:vv_vel} is implicit. 
Boris introduced a simple, geometrical procedure to find $\mathbf{v}_{n+1}$~\cite{boris1970}.
We use his trick in a slightly different way than usual, as a generic solver for an equation of the form
\begin{equation}
	\label{eq:boris_trick_form}
	\mathbf{v}_{n+1} = \mathbf{v}_n + \alpha \mathbf{E} + \beta \frac{\mathbf{v}_n + \mathbf{v}_{n+1}}{2} \times \mathbf{B} + \mathbf{c},
\end{equation}
where $\alpha$ and $\beta$ are some given scalar parameters and $\mathbf{E}$, $\mathbf{B}$ and $\mathbf{c}$ are some given vectors.
Note that $\mathbf{c}$ does not normally feature in most variants of the Boris integrator.
However, we will need it later as a ``container'' for various terms that arise from the Boris-SDC iteration.
Typically, $\alpha$ and $\beta$ are equal to the time step $\Delta t$ but we will generalise this when deriving the relativistic variant of Boris-SDC.
Furthermore, $\mathbf{E}$ will be the average of the electric fields at $\mathbf{x}_n$ and $\mathbf{x}_{n+1}$, $\mathbf{B}$ the magnetic field evaluated at some specific position whereas $\mathbf{c}$ will collect terms related, e.g., to the quadrature needed in Boris-SDC.
When used in this form, Boris' trick becomes Algorithm~\ref{alg:boris}.
Note that there a other slightly different variants.
In the terminology used by Zenitani and Umeda, we use the Boris-B algorithm~\cite{Zenitani_2018}.

\subsubsection{Interpolation of particle data to grid}
To compute $\mathbf{E}(\mathbf{x})$ in~\eqref{eq:lorentz}, we need to calculate the electrical field at the position of a particle from the mesh point values computed in Subsection~\ref{subsubsec:interpolation}.
To do so, linear interpolation is performed to collect the corresponding field value at each particle position. 
The value will be a sum of contributions from the surrounding nodes, each node contributing a field strength equal to the electric field at the node weighted by Eq.~\ref{eq:wLinear}. 
Any imposed, background electric and magnetic field can either be added to the nodes and interpolated or evaluated directly at the particle positions. For the simulations in this study, a neutralising static background electric field is imposed at the grid nodes. 
In the electrostatic case, this procedure is repeated in every time step.

\subsection{Boris-SDC\label{sec:3-prtcl-pushers}}
Boris-SDC is a time integration scheme for~\eqref{eq:lorentz} that provides tuneable order of accuracy. 
There are two slightly different versions. 
The one by Winkel et al.~\cite{winkel2015highOrderBoris} involves a substitution for velocity in the position update which can improve accuracy.
The substitution was dropped in the second variant by Tretiak and Ruprecht~\cite{TretiakRuprecht2019} to allow for the use of a GMRES-based convergence acceleration technique.
Both variants are based on collocation: the differential equation~\eqref{eq:lorentz} is turned into an integral equation
\begin{subequations}
	\begin{align}
		\mathbf{x}(t) &= \mathbf{x}_0 + \int_{t_0}^{t} \mathbf{v}(s)~ds, \\
		\mathbf{v}(t) &= \mathbf{v}_0 + \int_{t_0}^{t} \mathbf{f}(\mathbf{x}(s), \mathbf{v}(s))~ds,
	\end{align}
	\label{eq:lorentz_int}
\end{subequations}
where $\mathbf{x}_0$, $\mathbf{v}_0$ are approximations of $\mathbf{x}(t_0)$, $\mathbf{v}(t_0)$ brought forward from the previous time step.
In the formulation by Winkel et al.~\cite{winkel2015highOrderBoris}, the second equation is substituted into the first one so that
\begin{subequations}
\begin{align}
	\mathbf{x}(t) &= \mathbf{x}_0 + (t - t_0) \mathbf{v}_0 + \int_{t_0}^{t} \int_{t_0}^{r} \mathbf{f}(\mathbf{x}(s), \mathbf{v}(s))~ds~dr,  \\
	\mathbf{v}(t) &= \mathbf{v}_0 + \int_{t_0}^{t} \mathbf{f}(\mathbf{x}(s), \mathbf{v}(s))~ds.
\end{align}
\label{eq:lorentz_int_sub}
\end{subequations}
To compute an update from $t_0$ to $t_0 + \Delta t =: t_1$, the integrals are approximated using quadrature with respect to nodes $t_0 \leq \tau_1 < \ldots < \tau_M \leq t_1$.
Letting $\mathbf{x}_m$, $\mathbf{v}_m$ denote approximations for $\mathbf{x}(\tau_m)$, $\mathbf{v}(\tau_m)$ for $m = 1, \ldots, M$, these approximations read
\begin{subequations}
\begin{align}
	  \int_{t_0}^{t_1} \mathbf{v}(s)~ds &\approx \sum_{m=1}^{M} q_{m} \mathbf{v}_m,\\
	\int_{t_0}^{t_1} \mathbf{f}(\mathbf{x}(s), \mathbf{v}(s))~ds &\approx \sum_{m=1}^{M} q_{m} \mathbf{f}(\mathbf{x}_m, \mathbf{v}_m),
\end{align}
\end{subequations}
where the $q_m$ are quadrature weights.
Equations for the approximate values $\mathbf{x}_m$, $\mathbf{v}_m$ can be derived by inserting $t = \tau_m$ into~\eqref{eq:lorentz_int_sub} obtaining
\begin{subequations}
\begin{align}
	\mathbf{x}_m &= \mathbf{x}_0 + \mathbf{v}_0 \sum_{j=1}^{m} \Delta \tau_j +  \sum_{j=1}^{M} q_{m, j} \sum_{k=1}^{M} q_{j, k} \mathbf{f}(\mathbf{x}_k, \mathbf{v}_k), \\
	\mathbf{v}_m &= \mathbf{v}_0 + \sum_{j=1}^{M} q_{m, j} \mathbf{f}(\mathbf{x}_j, \mathbf{v}_j),
\end{align}
\end{subequations}
where $\Delta \tau_j := \tau_j - \tau_{j-1}$ for $j=1, \ldots, M$ and the $q_{m,j}$ are quadrature weights to approximate integrals $\int_{t_0}^{\tau_m} (\cdot)~ds$.
By subtracting the equations for index $m$ and $m-1$, the equations can be written in a node-to-node form
\begin{subequations}
\begin{align}
	\mathbf{x}_m &= \mathbf{x}_{m-1} + \Delta \tau_m \mathbf{v}_0 + \sum_{j=1}^{M} sq_{m,j} \mathbf{f}(\mathbf{x}_k, \mathbf{v}_k), \\
	\mathbf{v}_m &= \mathbf{v}_{m-1} + \sum_{j=1}^{M} s_{m, j} \mathbf{f}(\mathbf{x}_j, \mathbf{v}_j),
\end{align}
\label{eq:node-to-node-collocation}
\end{subequations}
where $s_{m, j} = q_{m, j} - q_{m-1,j}$ and the $sq_{m, j}$ can be calculated by rearranging the sum 
\begin{equation}
	\sum_{j=1}^{M} \left( q_{m, j} - q_{m-1,j} \right) \sum_{k=1}^{M} q_{j, k} = \sum_{j=1}^{M} s_{m, j} \sum_{k=1}^{M} q_{j, k},  
\end{equation}
see the Appendix in Winkel et al.~\cite{winkel2015highOrderBoris} for details.
Note that the equations for the $\mathbf{x}_m$, $\mathbf{v}_m$ are all coupled so that solving for them directly would require using a Newton iteration for a very large system of equations.
Instead, Boris-SDC computes approximations via a different iterative scheme where updates can be computed by a ``sweep'' of normal Boris integrator steps.

Skipping the derivation, which can also be found in Winkel et al.~\cite{winkel2015highOrderBoris}, the Boris-SDC iteration reads
\begin{subequations}
\begin{align}
\mathbf{x}_{m}^{k+1} &=\mathbf{x}_{m-1}^{k+1} +\Delta\tau_{m}\mathbf{v}_{0} 
   + \sum_{j=1}^{m-1} s_{m, j}^{X}\big (\mathbf{f} (\mathbf{x}_{j}^{k+1}, \mathbf{v}_{j}^{k+1}) -\mathbf{f} (\mathbf{x}_{j}^{k}, \mathbf{v}_{j}^{k}) \big) 
 +\sum_{j=1}^{M} sq_{m, j} \mathbf{f} (\mathbf{x}_{j}^{k}, \mathbf{v}_{j}^{k}) \label{eq:vvSDCx} \\
 \mathbf{v}_{m}^{k+1}&=\mathbf{v}_{m-1}^{k+1} + \frac{\Delta\tau_{m}}{2}\left( \mathbf{f} (\mathbf{x}_{m}^{k+1}, \mathbf{v}_{m}^{k+1}) - \mathbf{f} (\mathbf{x}_{m}^{k}, \mathbf{v}_{m}^{k}) \right),\label{eq:vvSDCv} \\
 &\quad + \frac{\Delta\tau_{m}}{2}\left( \mathbf{f} (\mathbf{x}_{m-1}^{k+1}, \mathbf{v}_{m-1}^{k+1}) - \mathbf{f} (\mathbf{x}_{m-1}^{k}, \mathbf{v}_{m-1}^{k}) \right) + \sum_{j=1}^{M} s_{m,j} \mathbf{f}(\mathbf{x}_j, \mathbf{v}_j), \ m=1, \ldots, M, \nonumber
\end{align}
\label{eq:vvSDC}
\end{subequations}
with $k$ counting iterations. 
The coefficients $s^X_{m, j}$ can be calculated from the distances $\Delta \tau_j$ between quadrature nodes, see again the Appendix in Winkel et al.~\cite{winkel2015highOrderBoris}.
If the iteration converges and $\mathbf{x}^{k+1}_m - \mathbf{x}^k_m \to 0$ and $\mathbf{v}^{k+1}_m - \mathbf{v}^k_m \to 0$, equations~\eqref{eq:vvSDC} reduce to the collocation equation~\eqref{eq:node-to-node-collocation}. 
Note that the position update is explicit: if we know the values from the previous iteration $k$ and all the values up to $\mathbf{x}^{k+1}_{m-1}$, we can directly compute $\mathbf{x}^{k+1}_m$ and so on.
In contrast, the velocity update is implicit, but we can use Boris' trick to compute $\mathbf{v}^{k+1}_m$.
To avoid cluttering the notation we assume that the charge-to-mass ratio $q/m$ is equal to unity here.
If that is not the case, just multiply the $\Delta \tau_m$ factors in front of the electric and magnetic field terms by $q/m$.
Let
\begin{equation}
\mathbf{c}^k_m :=- \frac{\Delta\tau_{m}}{2} \mathbf{f}(\mathbf{x}^k_m, \mathbf{v}^k_m) - \frac{\Delta \tau_m}{2}  \mathbf{f}(\mathbf{x}_{m-1}^k, \mathbf{v}_{m-1}^k) +\sum_{l=1}^{M} s_{m, l}\mathbf{f} (\mathbf{x}_{l}^{k}, \mathbf{v}_{l}^{k}) 
\end{equation}
collect all terms from the previous iteration with index $k$.
Expanding $\mathbf{f}$, the velocity update in the Boris-SDC iteration then reads
\begin{equation}
	\mathbf{v}^{k+1}_m = \mathbf{v}^{k+1}_{m-1} + \Delta \tau_m \frac{\mathbf{E}(\mathbf{x}_{m-1}^{k+1}) + \mathbf{E}(\mathbf{x}_m^{k+1})}{2} + \Delta \tau_m \frac{\mathbf{B}(\mathbf{x}_{m-1}^{k+1}) \times \mathbf{v}_{m-1}^{k+1} + \mathbf{B}(\mathbf{x}^{k+1}_m) \times \mathbf{v}_m^{k+1}}{2} + \mathbf{c}^k_m.
\end{equation}
To bring this into the form~\eqref{eq:boris_trick_form}, we add $-\mathbf{B}(\mathbf{x}_m^{k+1}) \times \mathbf{v}_{m-1}^{k+1} + \mathbf{B}(\mathbf{x}_m^{k+1}) \times \mathbf{v}_{m-1}^{k+1} = 0$ and let
\begin{equation}
	\tilde{\mathbf{c}}_m^k := \Delta \tau_m \frac{-\mathbf{B}(\mathbf{x}_m^{k+1}) \times \mathbf{v}_{m-1}^{k+1} + \mathbf{B}(\mathbf{x}_{m-1}^{k+1}) \times \mathbf{v}_{m-1}^{k+1}}{2} + \mathbf{c}_m^k
\end{equation}
so that, setting $\mathbf{E} := \frac{1}{2} \left( \mathbf{E}(\mathbf{x}^{k+1}_{m-1}) + \mathbf{E}(\mathbf{x}^{k+1}_m) \right)$ and $\mathbf{B} := \mathbf{B}(\mathbf{x}^{k+1}_m)$, the velocity update becomes
\begin{equation}
	\mathbf{v}_m^{k+1} = \mathbf{v}_{m-1}^{k+1} + \Delta \tau_m \mathbf{E} + \Delta \tau_m \frac{\mathbf{v}_{m-1}^{k+1} + \mathbf{v}_m^{k+1}}{2} \times \mathbf{B} + \tilde{\mathbf{c}}^k_m.
\end{equation}
This can now be solved using Algorithm~\ref{alg:boris}.
One time step of Boris-SDC then consists of the following steps:
\begin{enumerate}
	\itemsep0em
	\item Initialise $\mathbf{x}^0_m = \mathbf{x}_0$ and $\mathbf{v}^0_m = \mathbf{v}_0$ for $m=1, \ldots, M$.
	\item Perform $K$ sweeps:
		\begin{enumerate}
			\itemsep0em
			\item Evaluate $\mathbf{f}(\mathbf{x}^k, \mathbf{v}_m^k)$ for $m=1, \ldots, M$.
			\item Update $\mathbf{x}_m^{k} \rightarrow \mathbf{x}_{m}^{k+1}$ and $\mathbf{v}_m^k \rightarrow \mathbf{v}_m^{k+1}$ for $m=1, \ldots, M$ using~\eqref{eq:vvSDC}.
		\end{enumerate}
	\item If $\tau_M = t_{n+1}$, that is the end of the step is a quadrature node, set $\mathbf{x}^K_M \rightarrow \mathbf{x}_0$ and $\mathbf{v}_M^K \rightarrow \mathbf{v}_0$ and start the next time step.
\end{enumerate}
Throughout this paper, we use Gauss-Lobatto nodes for quadrature.
If other nodes are used, a final quadrature step is needed to deliver the approximate value at $t_{n+1}$.

\subsection{Relativistic Boris-SDC}\label{sec:3-rel-boris}
The relativistic Newton-Lorentz system in cgs units reads
\begin{subequations}
\begin{align}
\frac{d\mathbf{x}}{dt} &=\frac{\mathbf{u}}{\gamma(\mathbf{u})} =: \mathbf{g}(\mathbf{u}), \\
\frac{d\mathbf{u}}{dt} &= \frac{q}{m}\left(\mathbf{E}(\mathbf{x})+\frac{\mathbf{g}(\mathbf{u})}{c}\times\mathbf{B}(\mathbf{x})\right) =: \mathbf{f}(\mathbf{x}, \mathbf{u}),
\end{align}
    \label{eq:lorentz_rel}
\end{subequations}
where $\mathbf{u}$ is the proper velocity (the spatial component of the four-velocity) with the function $\mathbf{g}(\mathbf{u})$ yielding the coordinate velocity $\mathbf{v}$. 
Following Griffiths~\cite{griffiths2005introduction}, the relativistic Lorentz factor calculated from the proper velocity is
\begin{equation}
\gamma(\mathbf{u})=\sqrt{1+\mathbf{u} \cdot \mathbf{u}/c^{2}}.
\end{equation}
Written in this form, the system has the structure of a general second order initial value problem~\cite{HairerEtAl2003}. 
We consider the proper velocity $\mathbf{u}$ as the variable to solve for and, as far as the time stepping scheme is concerned, treat the coordinate velocity $\mathbf{v}$ as an auxiliary quantity.
In integral form,~\eqref{eq:lorentz_rel} becomes
\begin{subequations}
\begin{align}
	\mathbf{x}(t) &= \mathbf{x}_0 + \int_{t_0}^{t} \mathbf{g}(\mathbf{u}(s))~ds, \\
	\mathbf{u}(t) &= \mathbf{u}_0 + \int_{t_0}^{t} \mathbf{f}(\mathbf{x}(s), \mathbf{u}(s))~ds.
\end{align}
\label{eq:lorentz_rel_int}
\end{subequations}
Substituting as for the non-relativistic case would result in
\begin{equation}
    \mathbf{x}(t) = \mathbf{x}_0 + \int_{t_0}^{t} \mathbf{g}\left( \mathbf{u}_0 + \int_{t_0}^{r} \mathbf{u}(s)~ds \right) ~ dr.
\end{equation}
While this integral can be approximated by quadrature, it is not clear how the resulting SDC iteration can be written in a sweep-like fashion comparable to~\eqref{eq:vvSDC}.
We therefore use the less accurate formulation without substitution for the relativistic case and leave the derivation of a relativistic sweep with substitution for future work.
For the non-relativistic Newton-Lorentz equations~\eqref{eq:lorentz} the SDC sweep without substitution reads
\begin{equation}
\mathbf{x}_{m+1}^{k+1} = \mathbf{x}_{m}^{k+1}+\Delta\tau_{m}(\mathbf{v}_{m+1/2}^{k+1}-\mathbf{v}_{m+1/2}^{k}) + \sum_{j=1}^{M}s_{m,j}\mathbf{v}_{j}^{k},
\end{equation}
where 
\begin{equation}
\mathbf{v}_{m+1/2}^{k}=\mathbf{v}_{m}^{k}+\frac{\Delta\tau_{m}}{2}\mathbf{f}(\mathbf{x}_{m}^{k},\mathbf{v}_{m}^{k}).
\end{equation}
For the velocity, the iterations reads
\begin{equation}
\begin{aligned}
\mathbf{v}_{m+1}^{k+1} & =   \mathbf{v}_{m}^{k+1}+\frac{\Delta\tau_{m}}{2} \left[\mathbf{f}(\mathbf{x}_{m}^{k+1},\mathbf{v}_{m}^{k+1})+\mathbf{f}(\mathbf{x}_{m+1}^{k+1},\mathbf{v}_{m+1}^{k+1})\right]\\
 & -\frac{\Delta\tau_{m}}{2}\left[\mathbf{f}(\mathbf{x}_{m}^{k},\mathbf{v}_{m}^{k})+\mathbf{f}(\mathbf{x}_{m+1}^{k},\mathbf{v}_{m+1}^{k})\right] + \sum_{j=1}^{M}s_{m,j}\mathbf{f}(\mathbf{x}_{j}^{k},\mathbf{u}_{j}^{k}).
\end{aligned}
\end{equation}
For the relativistic Lorentz equations, the SDC sweep for the position becomes
\begin{equation}
\mathbf{x}_{m+1}^{k+1}=
\mathbf{x}_{m}^{k+1}+\Delta\tau_{m}\left[\mathbf{g}(\mathbf{u}_{m+1/2}^{k+1})-\mathbf{g}(\mathbf{u}_{m+1/2}^{k})\right] + \sum_{j=1}^{M}s_{m,j}\mathbf{g}(\mathbf{u}_{j}^{k}),
\end{equation}
where 
\begin{equation}
\mathbf{u}_{m+1/2}^{k+1}=\mathbf{u}_{m}^{k+1}+\frac{\Delta\tau_{m}}{2}\mathbf{f}\left(\mathbf{x}_{m}^{k+1},\mathbf{u}_{m}^{k+1}\right),
\end{equation}
and the sweep for the velocity
\begin{equation}
\begin{aligned}
\mathbf{u}_{m+1}^{k+1} & =  \mathbf{u}_{m}^{k+1}+\frac{\Delta\tau_{m}}{2}\left[\mathbf{f}(\mathbf{x}_{m}^{k+1},\mathbf{u}_{m}^{k+1})+\mathbf{f}(\mathbf{x}_{m+1}^{k+1},\mathbf{u}_{m+1}^{k+1})\right]\\
&  -\frac{\Delta\tau_{m}}{2}\left[\mathbf{f}(\mathbf{x}_{m}^{k},\mathbf{u}_{m}^{k})+\mathbf{f}(\mathbf{x}_{m+1}^{k},\mathbf{u}_{m+1}^{k})\right] + \sum_{j=1}^{M}s_{m,j}\mathbf{f}(\mathbf{x}_{j}^{k},\mathbf{u}_{j}^{k}).\\
\end{aligned}
\end{equation}
However, properly applying the Boris trick in the relativistic case requires some care.
Expanding $\mathbf{f}$ and defining
\begin{equation}
\mathbf{c}_{m}^k = \frac{\Delta\tau_{m}}{2c}\mathbf{g}(\mathbf{u}_{m}^{k+1})\times\mathbf{B}(\mathbf{x}_{m}^{k+1})-\frac{\Delta\tau_{m}}{2}\left[\mathbf{f}(\mathbf{x}_{m}^{k},\mathbf{u}_{m}^{k})+\mathbf{f}(\mathbf{x}_{m+1}^{k},\mathbf{u}_{m+1}^{k})\right] + \sum_{j=1}^{M}s_{m,j}\mathbf{f}(\mathbf{x}_{j}^{k}, \mathbf{u}_j^k)
\end{equation}
gives the update
\begin{equation}
\mathbf{u}_{m+1}^{k+1}=\mathbf{u}_{m}^{k+1}+\Delta\tau_{m} \left(\frac{\mathbf{E}(\mathbf{x}_{m}^{k+1})+\mathbf{E}(\mathbf{x}_{m+1}^{k+1})}{2}+\frac{1}{2c}\mathbf{g}(\mathbf{u}_{m+1}^{k+1})\times\mathbf{B}(\mathbf{x}_{m+1}^{k+1})\right)+\mathbf{c}_{m}^k.
\end{equation}
Setting
\begin{equation}
	\mathbf{E} := \frac{\mathbf{E}(\mathbf{x}_{m}^{k+1})+\mathbf{E}(\mathbf{x}_{m+1}^{k+1})}{2},
\end{equation}
and
\begin{equation}
	\mathbf{B} := \mathbf{B}(\mathbf{x}_{m+1}^{k+1}),
\end{equation}
results in
\begin{equation}
	\label{eq:rel_sdc_intermediate}
	\mathbf{u}_{m+1}^{k+1}=\mathbf{u}_{m}^{k+1}+\Delta\tau_{m}E+\frac{\Delta\tau_{m}}{2}\frac{\mathbf{u}_{m+1}^{k+1}}{c\gamma(\mathbf{u}_{m+1}^{k+1})}\times \mathbf{B}+\mathbf{c}_{m}^k.
\end{equation}
In the relativistic Boris algorithm, the Lorentz factor at the end of the time step must be estimated. 
The best way to estimate $\gamma(\mathbf{u}_{m+1})$ in the classical Boris schemes is still unclear~\cite{higuera2017structure,vay2008simulation}.
Typically, $\gamma$ is evaluated using the velocity after it has undergone half of the electric acceleration, that is $\gamma := \gamma(\mathbf{u}^{-})$.
This constant $\gamma$ is then included in the parameter $\beta$ in Algorithm~\ref{alg:boris}. 
Since $\mathbf{u}^{-}$ is not an approximation of the velocity at time $\tau_{m+1}$, when using this strategy in Boris-SDC, it will naturally not converge to $\mathbf{u}_{m}^{k+1}$ as $k$ increases.
Therefore, we would have
\begin{equation}
    \left\| \frac{\mathbf{u}^{k+1}_{m+1}}{c \gamma(\mathbf{u}^{-})} - \frac{\mathbf{u}^{k+1}_{m+1}}{c \gamma(\mathbf{u}^{k+1}_{m+1})} \right\| \geq \delta > 0
\end{equation}
for some $\delta$ that depends on the time step size and the size of the Lorentz factor for the problem but is independent of $k$.
This prevents the SDC iteration from converging to an accuracy smaller than $\delta$ because eventually, as $k$ increases, the constant error in the approximation of the relativistic factor will become dominant.

To prevent this and ensure convergence to the collocation solution we instead use $\gamma(\mathbf{u}_{m+1}^{k+1})\approx\gamma(\mathbf{u}_{m+1}^{k})$ to approximate $\gamma$.
This ensures that 
\begin{equation}
    \left\| \gamma(\mathbf{u}^{k+1}_{m+1}) - \gamma(\mathbf{u}^k_{m+1}) \right\| \to 0 \quad \text{as} \quad \left\| \mathbf{u}_{m+1}^{k}-\mathbf{u}_{m+1}^{k+1} \right\| \to 0
\end{equation}
so that the relativistic factor converges to its correct value as $k$ increases.
Letting $\gamma := \gamma(\mathbf{u}^{k}_m)$, Eq.~\eqref{eq:rel_sdc_intermediate} becomes
\begin{equation}
\mathbf{u}_{m+1}^{k+1}=\mathbf{u}_{m}^{k+1}+\Delta\tau_{m}E+\frac{\Delta\tau_{m}}{2\gamma c}\mathbf{u}_{m+1}^{k+1}\times \mathbf{B}+\frac{\Delta\tau_{m}}{2\gamma c}\left(\mathbf{u}_{m}^{k+1}\times \mathbf{B}-\mathbf{u}_{m}^{k+1}\times \mathbf{B} \right)+\mathbf{c}_{m}^k.
\end{equation}
Incorporating one term of the added zero into the constant by setting
\begin{equation}
	\tilde{\mathbf{c}}_{m}^k := \mathbf{c}_{m}^k - \frac{\Delta\tau_{m}}{2\gamma c}\mathbf{u}_{m}^{k+1}\times \mathbf{B},
\end{equation}
results in
\begin{equation}
	\mathbf{u}_{m+1}^{k+1}=\mathbf{u}_{m}^{k+1}+\Delta\tau_{m} \mathbf{E} + \frac{\Delta\tau_{m}}{\gamma c} \frac{\mathbf{u}_{m}^{k+1}+\mathbf{u}_{m+1}^{k+1}}{2}\times \mathbf{B} + \tilde{\mathbf{c}}_{m}^k.
\end{equation}
Now, the update step has again the correct form to be solved by Algorithm~\ref{alg:boris}.

%% file: results.tex
\section{Numerical Results\label{sec:4-results}}
We compare performance of Boris-SDC against the Boris integrator for two non-relativistic and two relativistic problems.
The first is a two-stream instability, representing a cold plasma with an exponentially growing instability.
The second is Landau damping, a hot plasma with an exponentially damped perturbation.
Both these problems are electrostatic and use a  one-dimensional PIC code.
Third, we compute a single relativistic particle in a Penning trap using an implementation of Boris-SDC in the Runko PIC software~\cite{runko}. The final experiment concerns a single relativistic particle in the special case where the magnetic and electric force exactly cancel out, where we compare the numerical drift for Boris-SDC against standard methods.

For the work-precision studies shown below, we compute the error between a given simulation and a reference simulation as the relative difference in the norm of the electric fields
\begin{equation}
	\mathrm{Err}^{\mathrm{rel}} = \frac{ | \left\| E_{\text{ref}} \right\|_{l_2} - \left\| E \right\|_{l_2} |}{ \left\| E_{ref} \right\|_{l_2} }
	\label{eq:error}
\end{equation}
where
\begin{equation}
||E||_{l_2}=\sqrt{\Delta z\sum_{i=1}^{N_{i}}E_{i}^{2}}, \label{eq:el2}
\end{equation}
is used to calculate the norm.
Here, $\Delta z$ is the grid spacing, $N_{i}$ is the number of cell nodes and $E_{i}$ the electric field at node $i$.
Note that we compute the relative error of the $l_2$-norms of the electric fields and not the relative error in the electric field directly.
This makes comparison of results on different mesh resolutions easier.

\subsection{Two-Stream Instability}
\begin{figure}[t]
    \centering
    \includegraphics[width=1\linewidth]{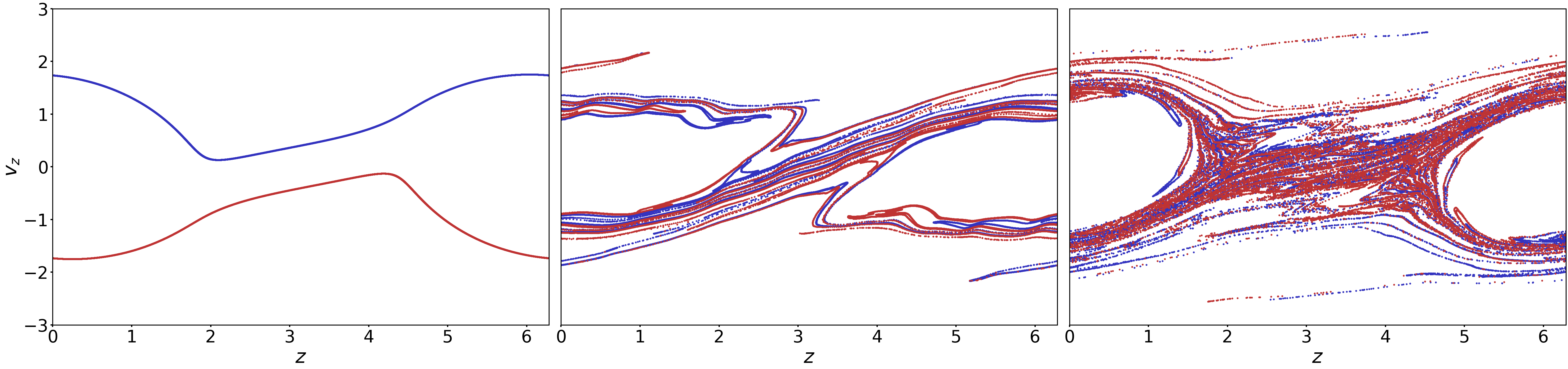}%
    \caption{Two-Stream Instability particle position-velocity phase-space at $t=60, \, 180, \, 300$
    for a half period sinusoidal perturbation in the charge density of
    each beam (magnitude $A=10^{-1}$).}
    \label{fig:tsi_dynamics}
\end{figure}

The two-stream instability is a type of streaming instability.
A beam containing one species of charged particles streams through another. 
While the instability can occur for counter-streaming beams of any particle mass and charge, for the purposes of this study we focus on beams consisting of the same species with particles of identical charge and mass.
Such counter-streaming beams are inherently unstable as any perturbation in density or velocity distribution is reinforced by the charge induced in the other beam and vice-versa~\cite{chen1974}. 
The dynamics become increasingly chaotic as the instability develops and the plasma becomes increasingly thermalised until a phase-space structure resembling an eye appears, see Figure~\ref{fig:tsi_dynamics}.
We start with a sinusoidally perturbed particle density distribution 
\begin{equation}
    n(x) = n_0 + A \cos \left( \frac{2\pi k x}{L} \right),
\end{equation}
where $n_0$ is the uniform distribution, $A$ is the initial perturbation magnitude, $k$ is the perturbation mode number and $L$ is the domain length.
Assuming that $A$ is small, $A \ll 1$, the early dynamics are linear with the electric field strength growing exponentially at rate \footnote{Note that this $\gamma$ is not the relativistic Lorentz factor.} $\gamma$.
This rate can be calculated analytically.
The dispersion relation of the electrostatic wave induced by the perturbation has four roots~\cite{birdsallLangdon1985} and the growth rate of the instability corresponds to the imaginary root
\begin{equation}
	\gamma_{\mathrm{theory}}=\left[k^{2}v_{0}^{2}+\omega_{p}^{2}-\omega_{p} (4k^{2}v_{0}^{2}+\omega_{p}^{2}) ^{1/2}\right]^{1/2}.
\end{equation}
Here, $k$ is the wave number of the perturbation, $v_{0}$ is the initial velocity magnitude of the beams and $\omega_{p}$ is the plasma frequency
\begin{equation}
	\omega_{p}=\sqrt{\frac{nq{{}^2}}{m\epsilon}}, 
\end{equation}
using the plasma density $n$, particle charge $q$ and mass $m$ as well as the permittivity $\epsilon$.

To show that Boris-SDC captures the early dynamics correctly, we simulate the initial growth of the electric field for a weak ($A=10^{-4}$) and strong ($A=10^{-1}$) perturbation. 
In both cases we expect the field to grow exponentially with rate $\gamma$, up to a saturation point, followed by transition to a chaotic, nonlinear regime.
The simulation uses $N_{q}=10^{4}$ particles, $N_{z}=100$ grid nodes and a time-step size $\Delta t=0.1$. 
A periodic domain of length $L=2\pi$ is used, a perturbation with wave number $k=1$, beam velocity is set to $|v_{0}|=1$, the particle/mass ratio  to $q/m=1$, and the permittivity to $\epsilon=1$. 
Particle charge was calculated so that $\omega_{p}=1$ using a plasma density defined by $n=N_{q}/L$.

\begin{figure}[t]
\begin{minipage}[c]{0.475\textwidth}%
\includegraphics[width=1\linewidth]{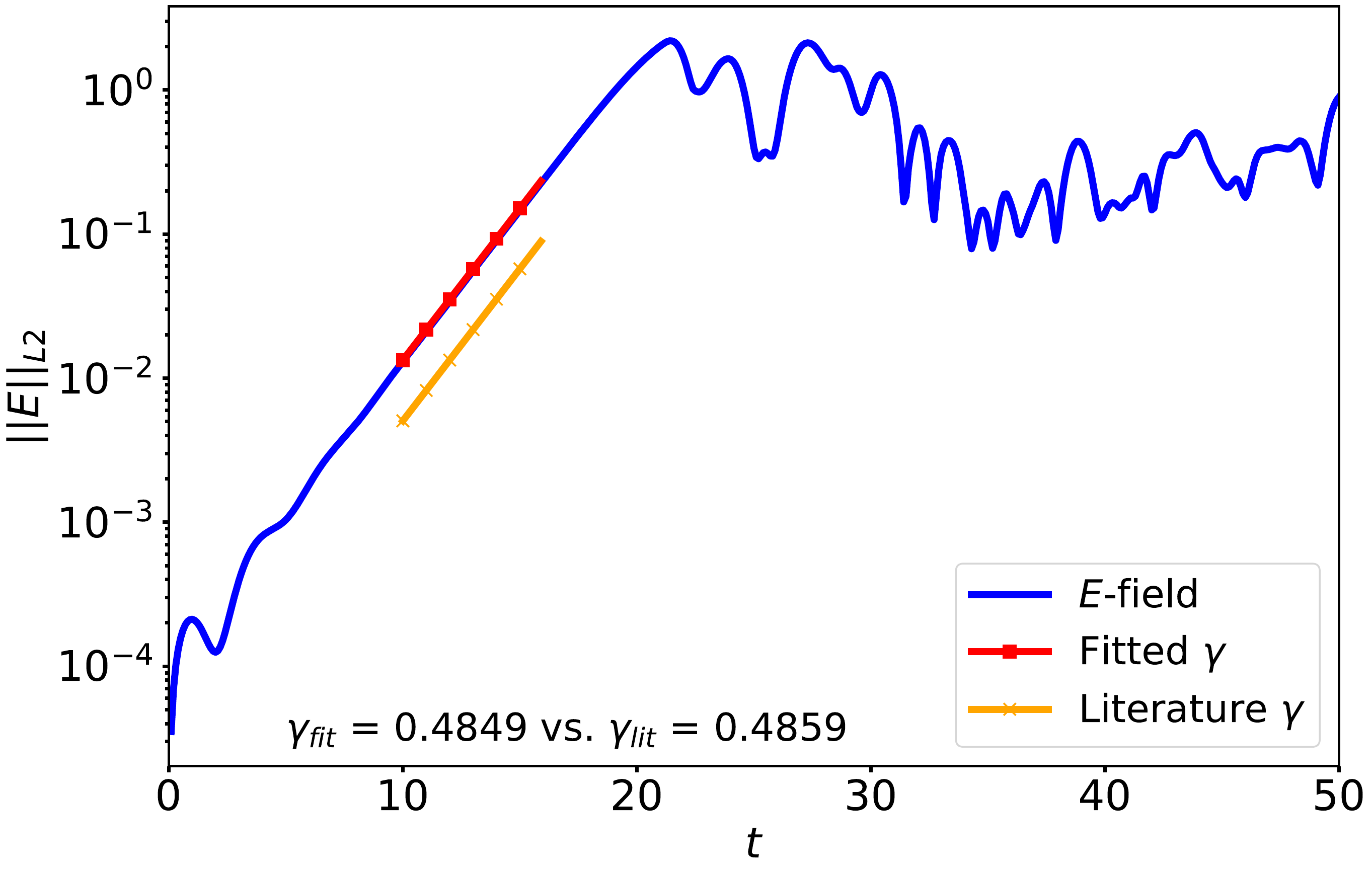}
\end{minipage}%
\begin{minipage}[c]{0.475\textwidth}%
\includegraphics[width=1\linewidth]{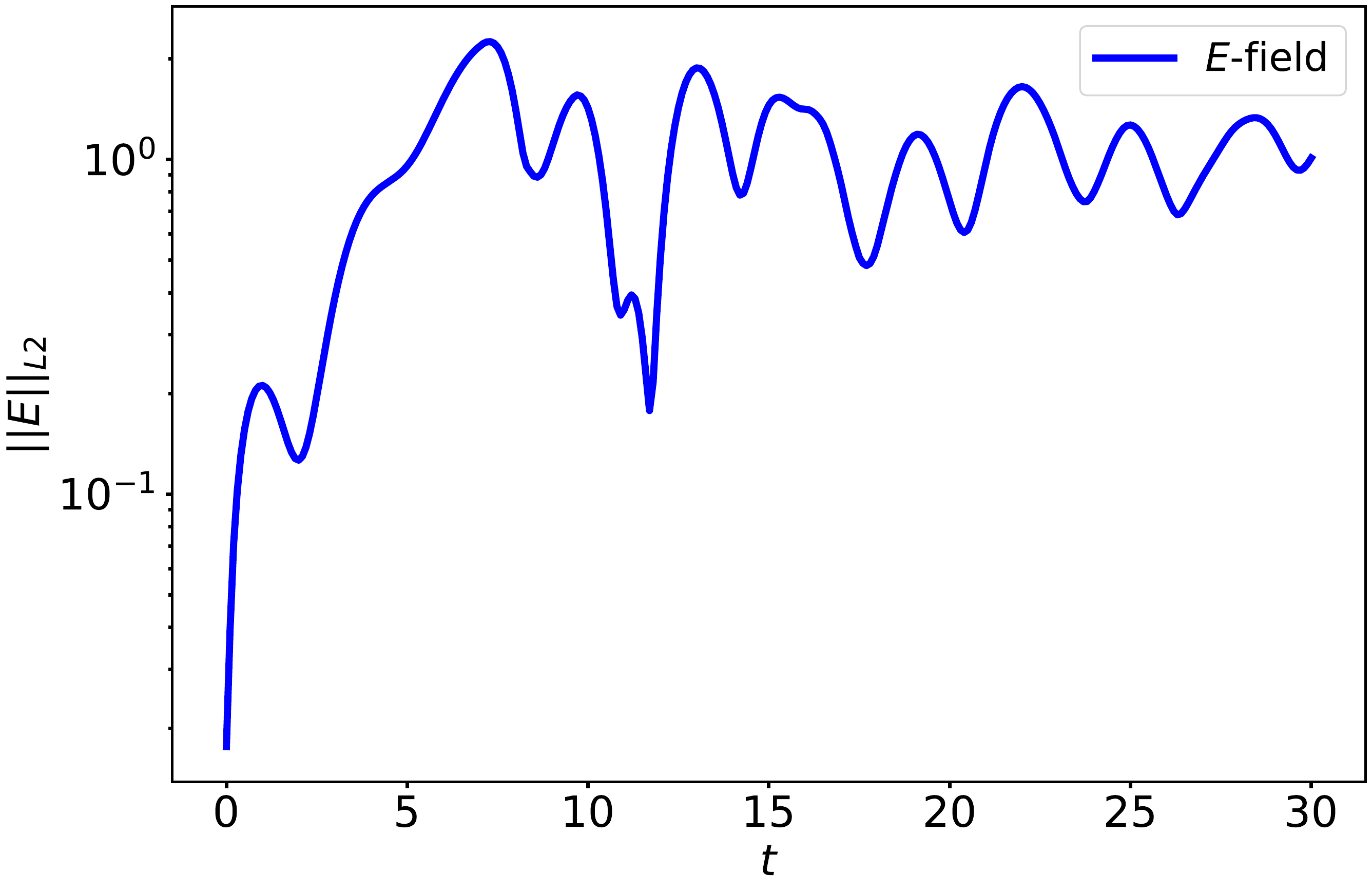}
\end{minipage}
\caption{$E$-field $l_2$ norm growth of weak (left) and strong (right) two-stream instability.}
\label{fig:tsi_growth}
\end{figure}

Figure~\ref{fig:tsi_growth} shows the electric field strength given by~\eqref{eq:el2} over time for the weak (left) and strong perturbation (right).
For the weak perturbation case, linear growth can be seen up to a time of around $t = \SI{20}{\second}$.
Between $t=\SI{10}{\second}$ and $t=\SI{16}{\second}$, a line with slope $\gamma_{\mathrm{theory}}$ is shown.
A best fit of the electric field data between \SI{12}{\second} and \SI{18}{\second} gives a growth rate of $\gamma_{\mathrm{fit}}=0.4849$ which matches the rate predicted by theory ($\gamma_{\mathrm{lit}}=0.4859$) to within $0.19\%$.
A simulation using classic Boris produces $\gamma_{\text{boris}} = 0.4855$.

To investigate performance we compare the computational effort in terms of right hand side evaluations required by Boris and Boris-SDC to reach a certain error.
To minimise noise so that we can clearly assess the error from numerical discretisation, we use $N_q = 2 \times 10^{5}$ particles in all simulations.
Meshes with $N_z = 10$, $N_z = 100$ and $N_z = 1000$ points were used to analyse the interplay of spatial and temporal discretisation errors.
Simulation parameters are summarised in Table~\ref{tab:tsi_wp_params}. 
The reference solution for both cases is a high accuracy Boris-SDC ($M=3$, $K=3$) simulation using 5 times the maximum time- and space resolution, thus $N_t=5000$ and $N_z=5000$.
Boris-SDC uses $M=3$ Gauss-Lobatto nodes so that the underlying quadrature is fourth order accurate.
\begin{table}[bht]
\centering{}\caption{Parameter used for work-precision study.}
\label{tab:tsi_wp_params} %
\begin{tabular}{cccc}
\hline 
\textbf{Parameter} & \textbf{Key} & \textbf{Values}\tabularnewline
\hline
Integrator & - & {[}Boris-SDC ($M=3$, $K=2$) , Boris{]} \tabularnewline
Time steps & $N_t$ & $[10, 20, 40, 50, 80, 100, 200, 400, 500, 1000]$ \tabularnewline
Mesh resolution & $N_z$ & $[10, 100, 1000]$ \tabularnewline
Particle count & $N_q$ & $2\cdot10^{5}$ \tabularnewline
End time & $T_{E}$ & $10$ \tabularnewline
\hline
\end{tabular}
\end{table}

\begin{figure}[t]
	\centering
	\includegraphics[width=0.475\linewidth]{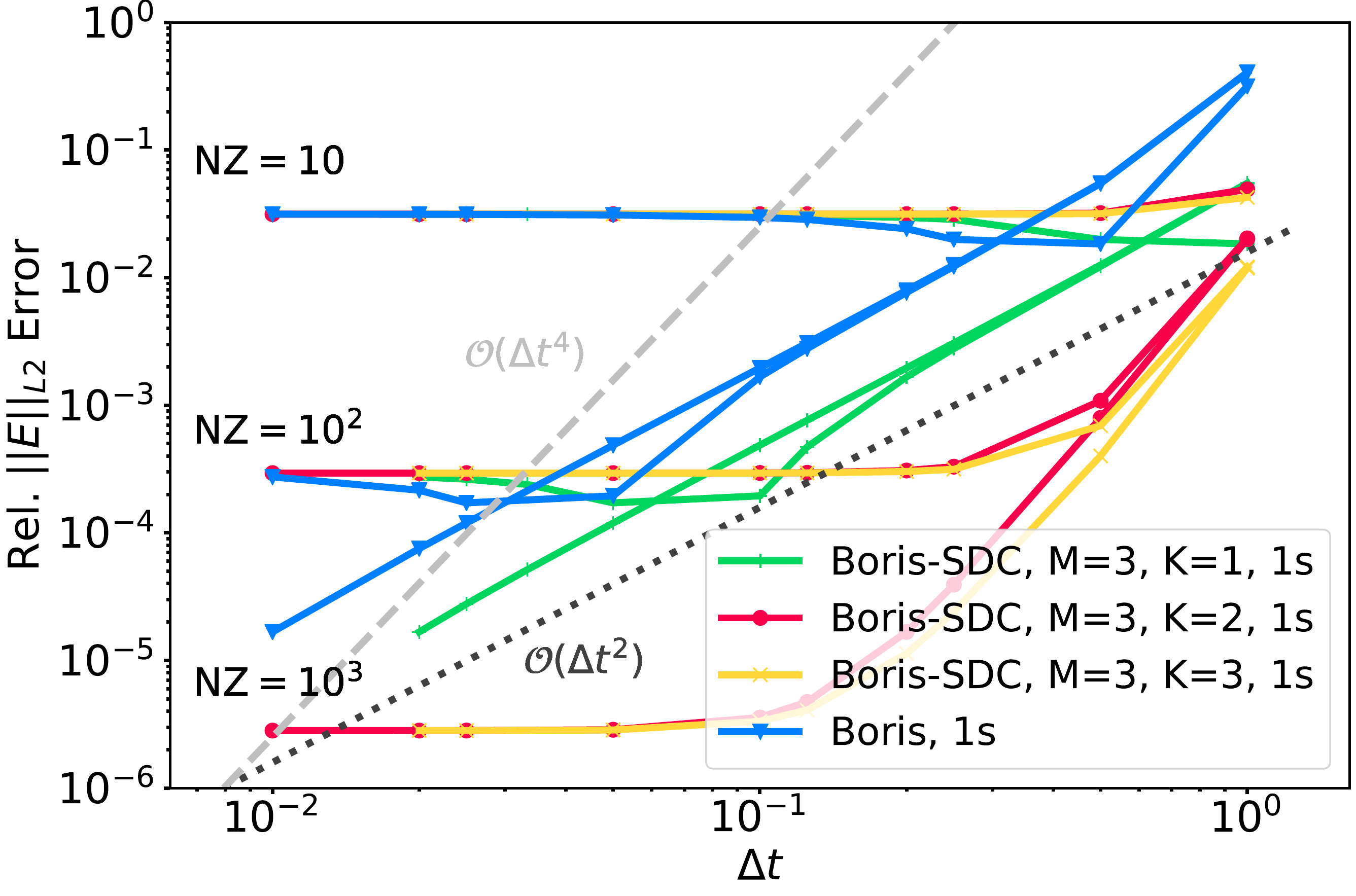}
	\includegraphics[width=0.475\linewidth]{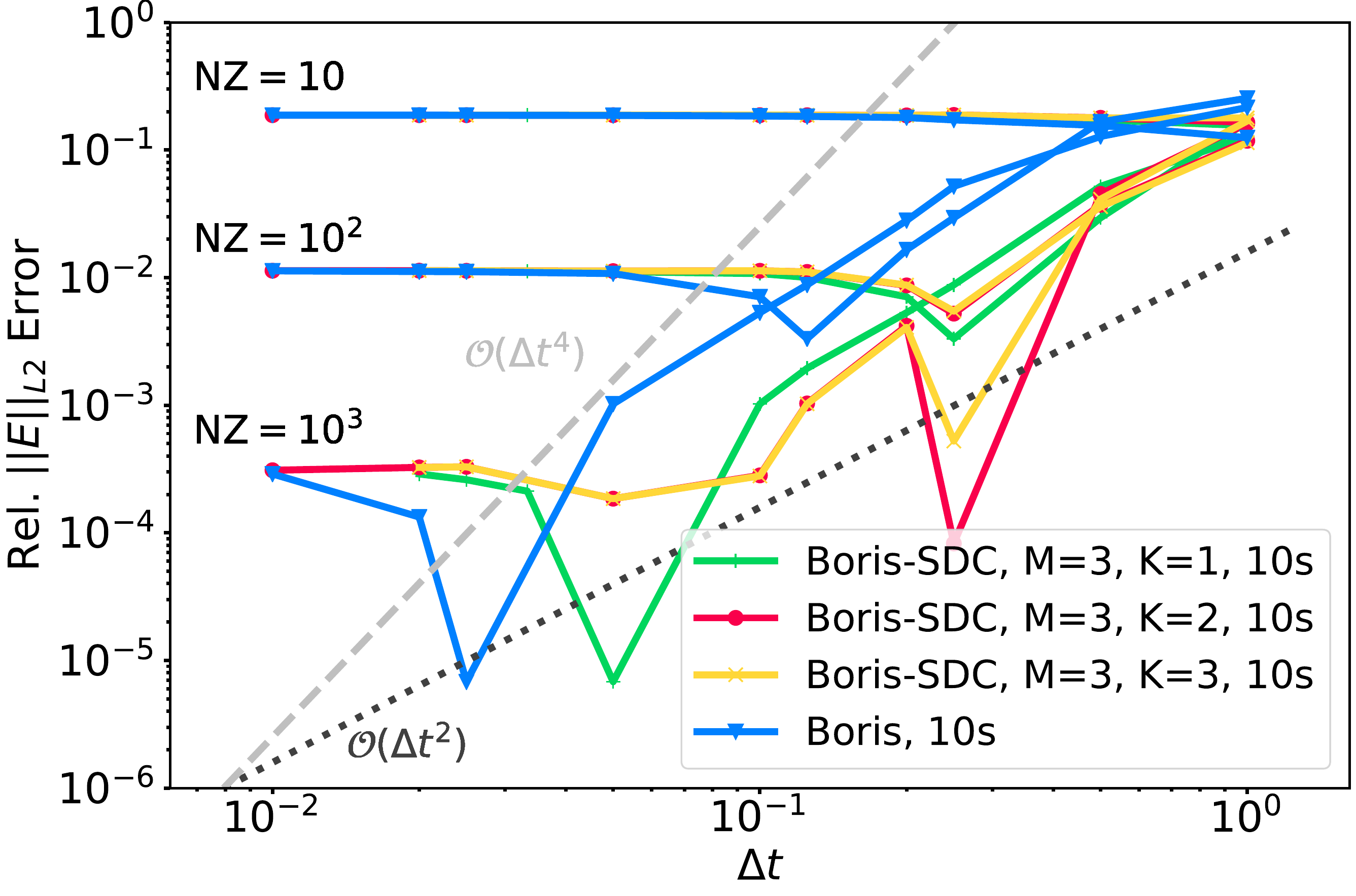}
	\caption{Accuracy comparison of Boris and Boris-SDC for the strongly perturbed two-stream instability at 1 and 10s simulation time across 3 Boris-SDC iteration counts. Dashed/dotted lines are guide lines for fourth and second order convergence respectively.}
\label{fig:tsi_strong_full_dt}
\end{figure}

Figure \ref{fig:tsi_strong_full_dt}~shows the relative error at \SI{1}{\second} (left) and \SI{10}{\second} (right) versus time-step size for the strongly perturbed two-stream instability.
At \SI{1}{\second} simulation time, both Boris-SDC with $K=1$ iteration and standard Boris are second order accurate, but Boris-SDC has a slightly smaller error constant.
For $K=2$ and $K=3$ iterations, Boris-SDC is fourth order accurate.
Here, the limiting factor is the order $p = 2 M - 2 = 4$ of the underlying quadrature rule.
By contrast, at \SI{10}{\second} simulation time, the dynamics are strongly nonlinear and convergence orders are less clear.
Accuracy is mostly determined by the spatial error: for $N_z = 10$ and $N_z = 100$ mesh points, there is little impact from the time step size for both Boris and Boris-SDC anymore.

\begin{figure}[t]
	\centering
	\includegraphics[width=0.475\linewidth]{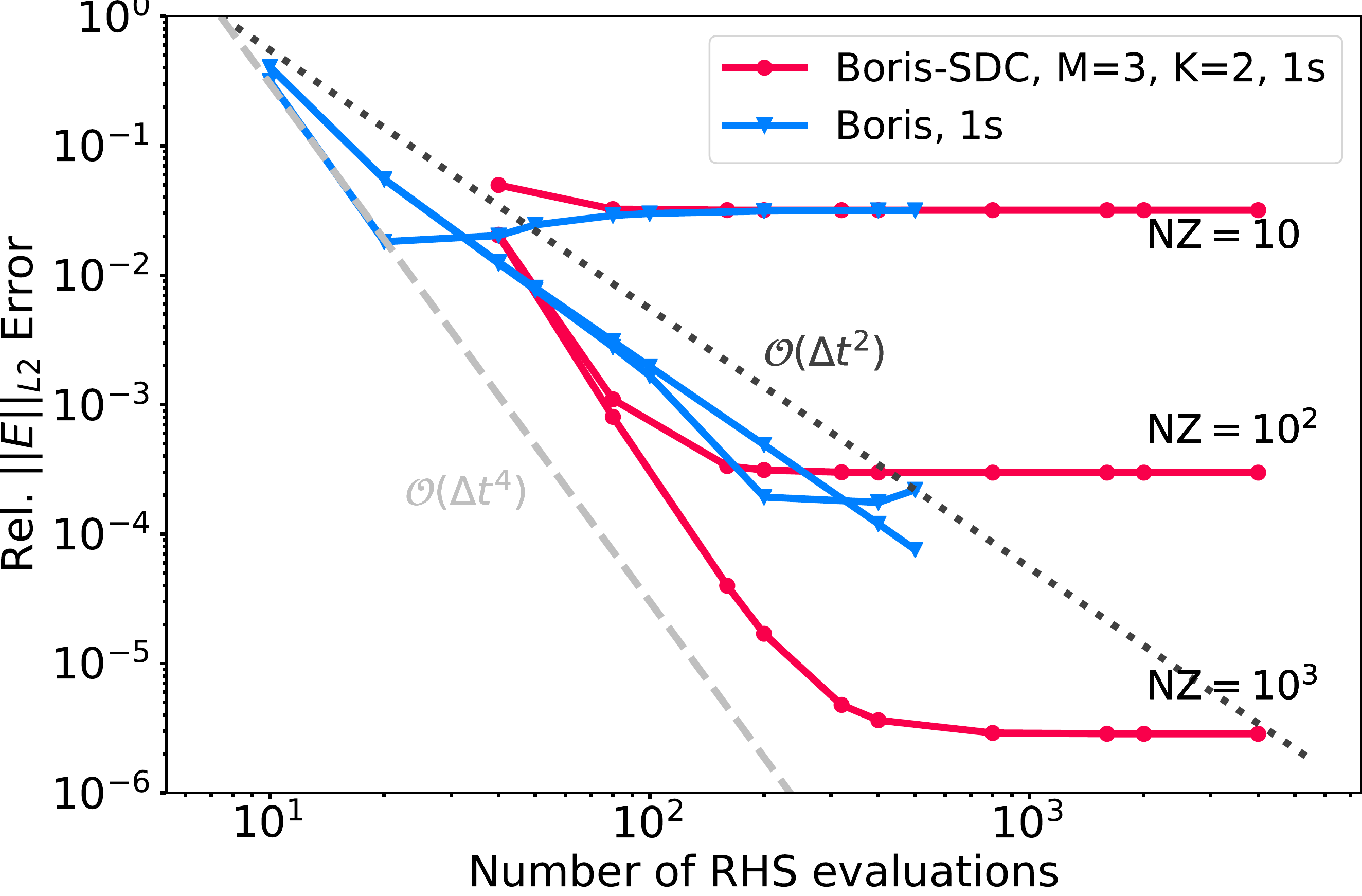}
	\includegraphics[width=0.475\linewidth]{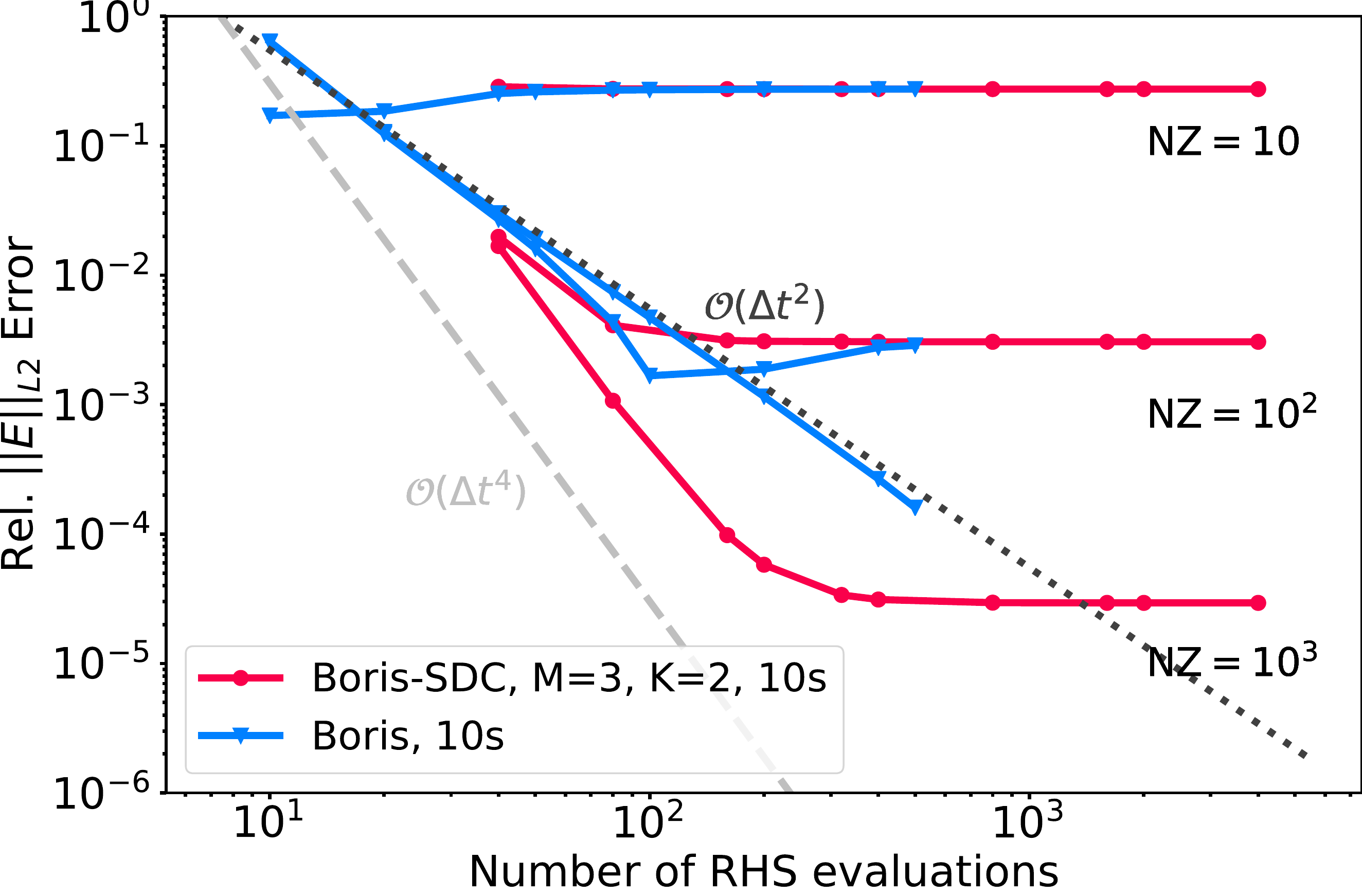}
	\caption{Performance comparison of Boris and Boris-SDC for the weakly perturbed two-stream 
instability at 1 and 10s simulation time. Dashed/dotted lines are
guide lines for fourth and second order convergence respectively.}
	\label{fig:tsi_weak_versus_rhs}
\end{figure}

Figure \ref{fig:tsi_weak_versus_rhs} shows error versus computational work, measured by the number of required right hand side evaluations, for the weakly perturbed two-stream instability at \SI{1}{\second} (left) and \SI{10}{second} (right).
For errors above 1\%, the Boris integrator is the more efficient choice as Boris-SDC will require more computational work.
If errors of 1\% or below are required, Boris-SDC becomes more efficient.
To achieve an error of, say, $10^{-4}$ at \SI{10}{\second} simulation time, Boris-SDC requires about 200 evaluations of the right hand side whereas Boris requires around 800.
Because spatial resolution remains fixed, eventually there are no more gains for both methods from decreasing $\Delta t$ further, as the error becomes dominated by the contributions from spatial errors.

\begin{figure}[t]
	\centering
	\includegraphics[width=0.475\linewidth]{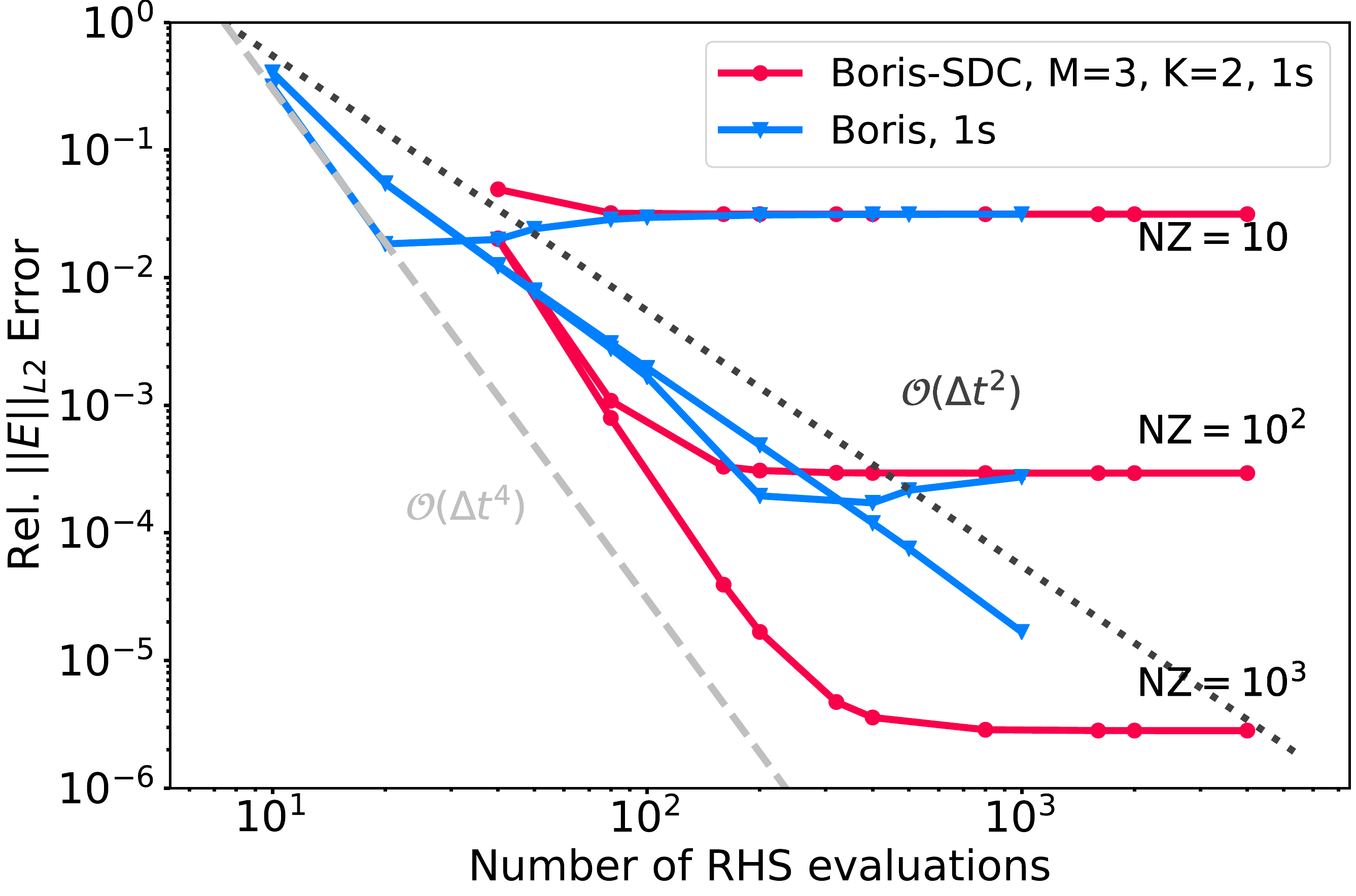}
	\includegraphics[width=0.475\linewidth]{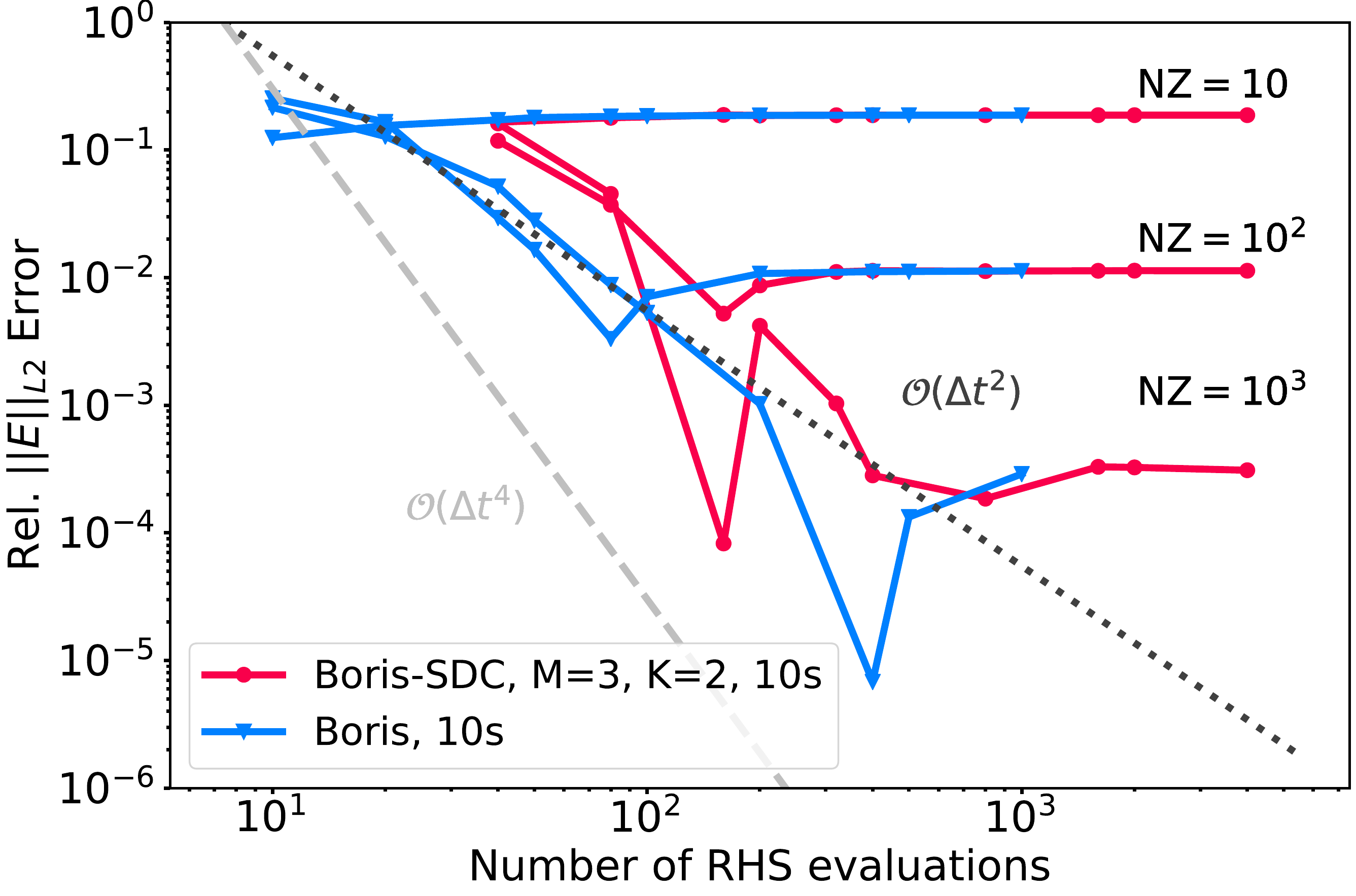}
	\caption{Performance comparison of Boris and Boris-SDC for the strongly perturbed two-stream
instability at 1 and 10s simulation time. Dashed/dotted lines are
guide lines for fourth and second order convergence respectively. }
	\label{fig:tsi_strong_versus_rhs}
\end{figure}

Figure~\ref{fig:tsi_strong_versus_rhs} shows error at \SI{1}{\second} (left) and \SI{10}{\second} (right) versus computational effort for the strongly perturbed case. 
At 1 second, performance is similar to the weakly perturbed case with Boris-SDC being somewhat more efficient than Boris.
At 10 seconds, however, the two integrators deliver comparable performance.
The higher order of Boris-SDC allows one to take fewer, larger time steps but, in contrast to the weakly perturbed case, this gain is off-set by the increased per-time step cost of Boris-SDC.

\subsection{Landau Damping}
\begin{figure}[t]
    \centering
    \includegraphics[width=1\linewidth]{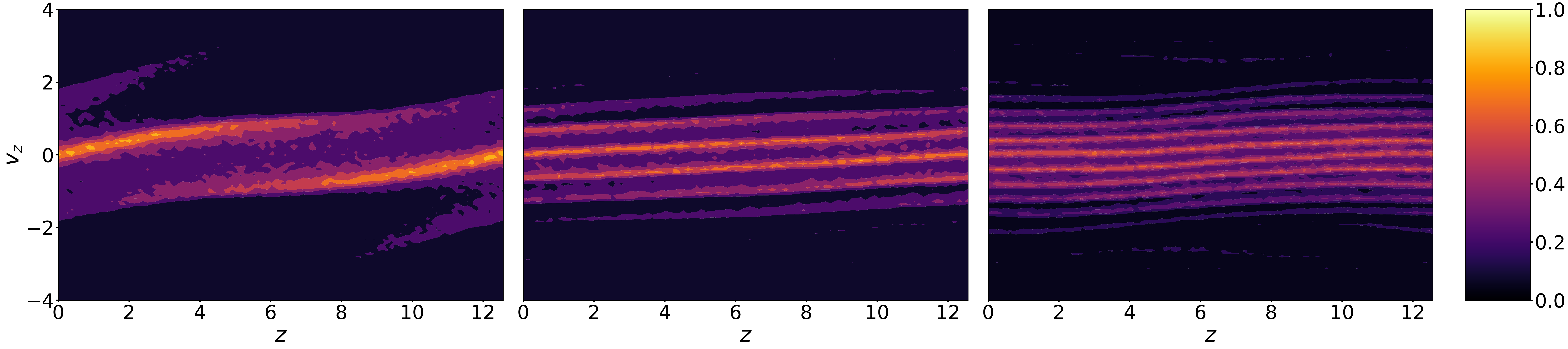}%
    \caption{Landau damping density distribution function at $t=60, \, 180, \, 300$ for a single period sinusoidal perturbation in the charge density of each beam (magnitude $A=0.5$).}
    \label{fig:landau_dynamics}
\end{figure}

Landau damping refers to the attenuation of electrostatic waves in a collisionless plasma from energy transfer between particles and the electric field. 
A physical description of the phenomenon is given by Chen~\cite{chen1974} while a detailed mathematical analysis can be found in the original study by Landau~\cite{landau1946}.
Figure~\ref{fig:landau_dynamics} shows the evolution of phase space density for a strongly perturbed density wave in a negatively charged particle distribution on a neutralising background. 
The plasma slowly returns towards an equilibrium state due to Landau damping.

As for the two-stream instability we first demonstrate that the code captures correctly the evolution of the electric field strength in both weakly and strongly perturbed simulations. 
Similar to the two-stream instability, the magnitude of an imposed sinusoidal density perturbation will determine the linearity of the ensuing wave damping. Both linear and non-linear Landau damping has been studied extensively~\cite{ayuso2012landau,canosa1972landau,cheng1976landau,nakamura1999landau,rossmanith2011landau} for initial density distributions $n, \, f$ in position $x$ and velocity $v$ space of the form
\begin{equation}
	\label{eq:landau_dens}
	n (x, 0) =\rho (x, 0) =1+A\cos (kx),
\end{equation}
and
\begin{equation}\label{eq:landau_vel}
f (v, 0) =\frac{1}{\sqrt{2\pi}}e^{-\frac{v^{2}}{2}}, 
\end{equation}
with magnitudes $A=0.01$ and $A=0.5$ for the weakly and strongly perturbed regimes respectively used in all studies. Note that alternatively the velocity could be perturbed as $v = v(x)$ instead of the position.

The studies mentioned above apply numerical methods directly to the Maxwell-Vlasov system, which allows direct use of~\eqref{eq:landau_dens} and~\eqref{eq:landau_vel} as initial conditions. 
To realise the same setup in an ESPIC code requires a particle distribution that corresponds to the initial density distribution. 
The average unit density of particles $n_{0}$ can be calculated from~\eqref{eq:landau_dens} using
\begin{equation}
	n_{0}=\frac{1}{L} \int^{L}n(x,0) \, dx = \frac{1}{L} \int^{L}1+A\cos(kx) \, dx.
\end{equation}
Since the charge density of the particle species is given globally, the charge of the macro-particles must be assigned based on this.
Unlike the two-stream instability, this problem does not use the plasma frequency as an independent quantity. 
The particle charge was calculated for a given simulation by dividing the global charge of the species by the desired number of computational particles~$N_q$:
\begin{equation}
	q=\frac{Q}{N_q}=\frac{\int^{L}\rho(x)\,dx}{N_q}=\frac{\int^{L}1+A\cos(kx)\,dx}{N_q}.
\end{equation}
Finally, the particle velocities can be distributed randomly to fit the Maxwellian defined by~\eqref{eq:landau_vel}. 

A perturbation mode $k=0.5$ and domain length $L=4\pi$ are used and leads to a single perturbation period in space. 
The expected damping rate from linear theory is $\gamma_{\mathrm{theory}}=-0.1533$, while the above studies report damping rates in the range $\gamma_{\mathrm{lit}}=[-0.292, -0.220]$ for the strongly perturbed setup.
Rapid oscillation of the electric field is expected with an overall exponential damping of the perturbation until a saturation point is reached. 
In the weakly perturbed case, the plasma should continue irregular oscillation after the saturation point, while a phase of slight growth should be observed in the strongly perturbed dynamics.
Simulation parameter are summarised in Table~\ref{tab:landau_params}.

\begin{table}[h]
    \centering
    \caption{Physical setup parameters for the Landau damping cases.}
    \label{tab:landau_params} %
    \begin{tabular}{ccc}
    \toprule 
    \textbf{Parameter} & \textbf{Weak perturbation} & \textbf{Strong perturbation}\tabularnewline
    \midrule 
    $L$ & $4\pi$ & $4\pi$\tabularnewline
    $k$ & $0.5$ & $0.5$\tabularnewline
    $A$ & $0.05$ & $0.5$\tabularnewline
    $v_{th}$ & $1$ & $1$\tabularnewline
    $\epsilon$ & $1$ & $1$\tabularnewline
    \midrule 
    $\omega_{p}$ & $1$ & $1$\tabularnewline
    $n_{0}$ & $1$ & $1$\tabularnewline
    $q$ & $L/N_q$ & $L/N_q$\tabularnewline
    \bottomrule
    \end{tabular}
\end{table}

Weak Landau damping is difficult to capture with PIC unless a large number of particles is used, since the driving interaction is the energy exchange between wave and particles in the trapping range, close to the phase velocity $v_{\phi}$. 
For the studied setup, the phase velocity is placed toward the tail end of the velocity distribution, meaning only a small proportion of the particles are in the trapping range. 
We found that using $A=0.01$ did not induce noticeable damping and thus we increased the magnitude of the perturbation to $0.05$ to place more particles in the trapping range.
Each simulation used $N_{q}=10^{5}$ particles, $N_{z}=100$ grid nodes and a time-step size $\Delta t=0.1$.

\begin{figure}[t]
    \centering
    \includegraphics[width=.49\linewidth]{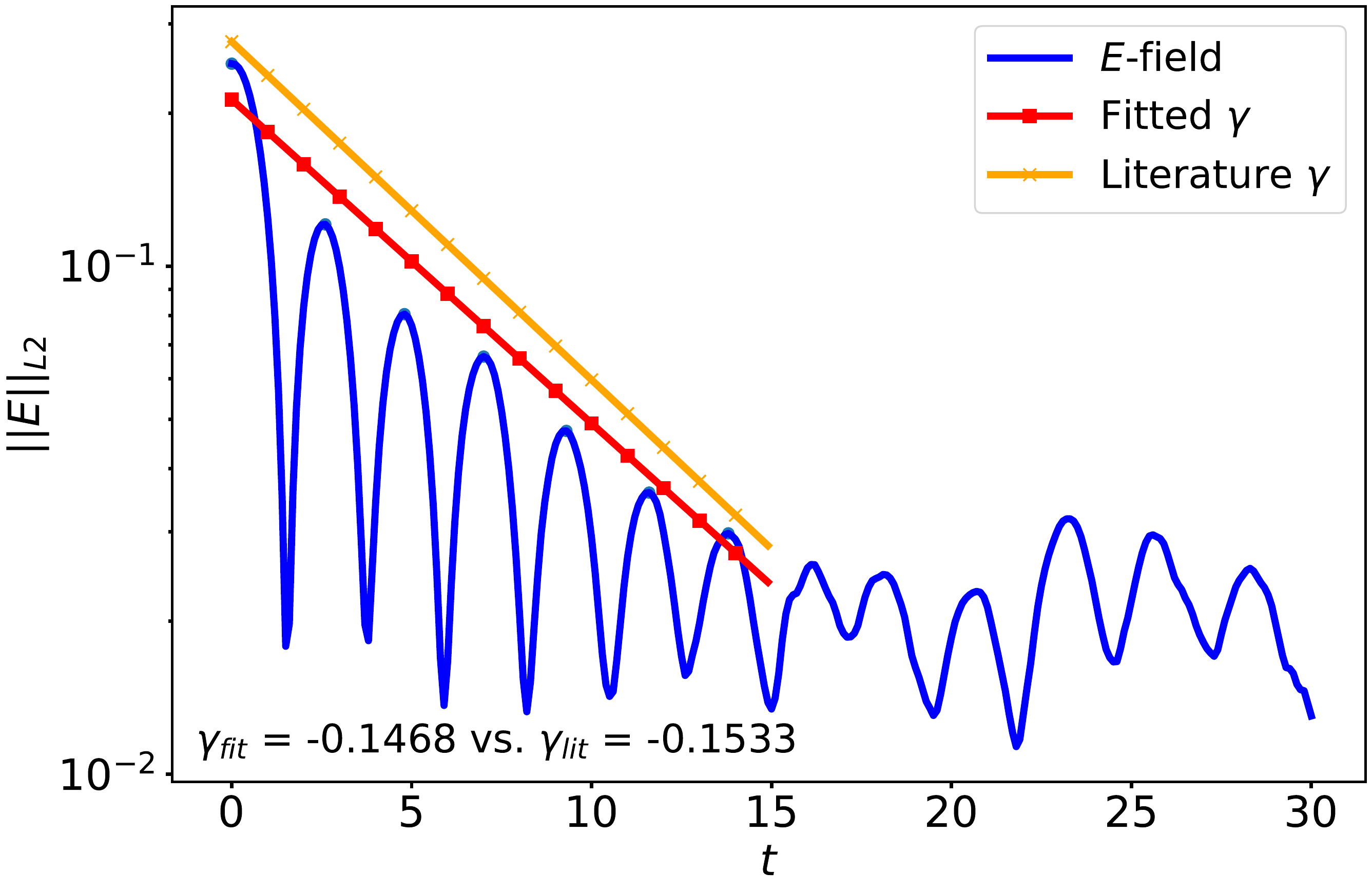}
    \includegraphics[width=.49\linewidth]{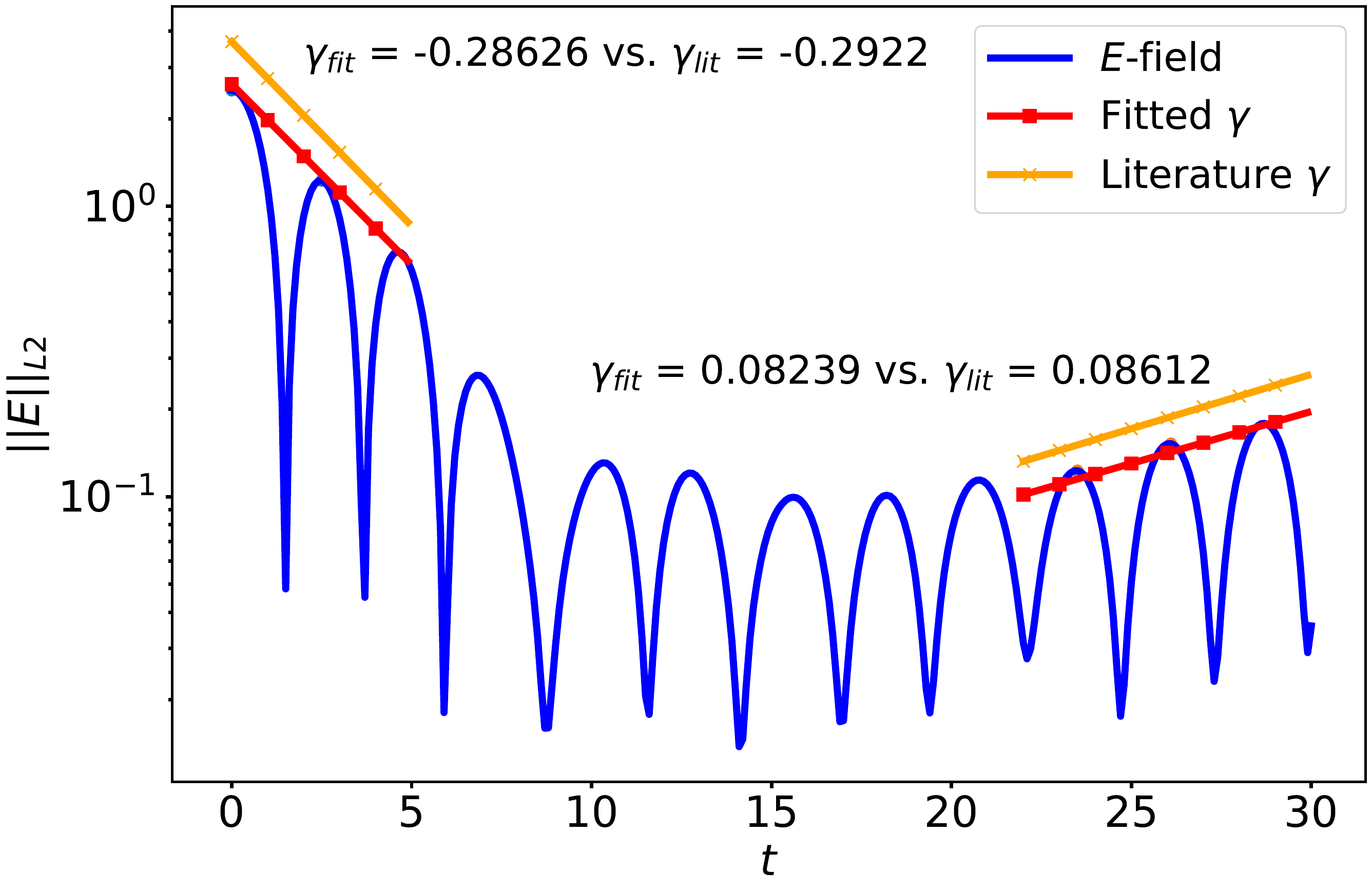}
    \caption{E-field $l_2$-norm evolution for weak (left) and strong (right) Landau damping.}
    \label{fig:landau_growth}
\end{figure}

Figure~\ref{fig:landau_growth} shows the evolution of the magnetic field for the weak (left) and strong (right) perturbation. 
The dynamics match the expected behaviour: for the weak perturbation, the electric field decreases exponentially until a saturation point at around $t = 15$.
After that, the field continues to oscillate without a clear change in magnitude.
For the strong perturbation, the damping phase is shorter and after saturation, a slight growth of the electric field sets in.

To evaluate the damping rates, best fit lines were drawn through the relevant oscillation peaks on each graph, the first seven peaks were used for the linear data ($t\sim[0, 15]$) and first three used for the non-linear data ($t\sim[0, 5]$). The relative error of the simulated linear damping rate $\gamma_{\mathrm{fit}}=-0.147$ to the theoretical value was found to be approximately $4. 1\%$, thought to be reasonable agreement considering the rapidly oscillating dynamics and difficulties related to the scheme. Furthermore, the simulated damping rate in the non-linear case $\gamma_{\mathrm{fit}}=-0.28825$ was firmly within the range of values reported in literature.

\begin{figure}[t]
    \centering
    \begin{minipage}{0.475\textwidth}
        \centering
        \includegraphics[width=1\linewidth]{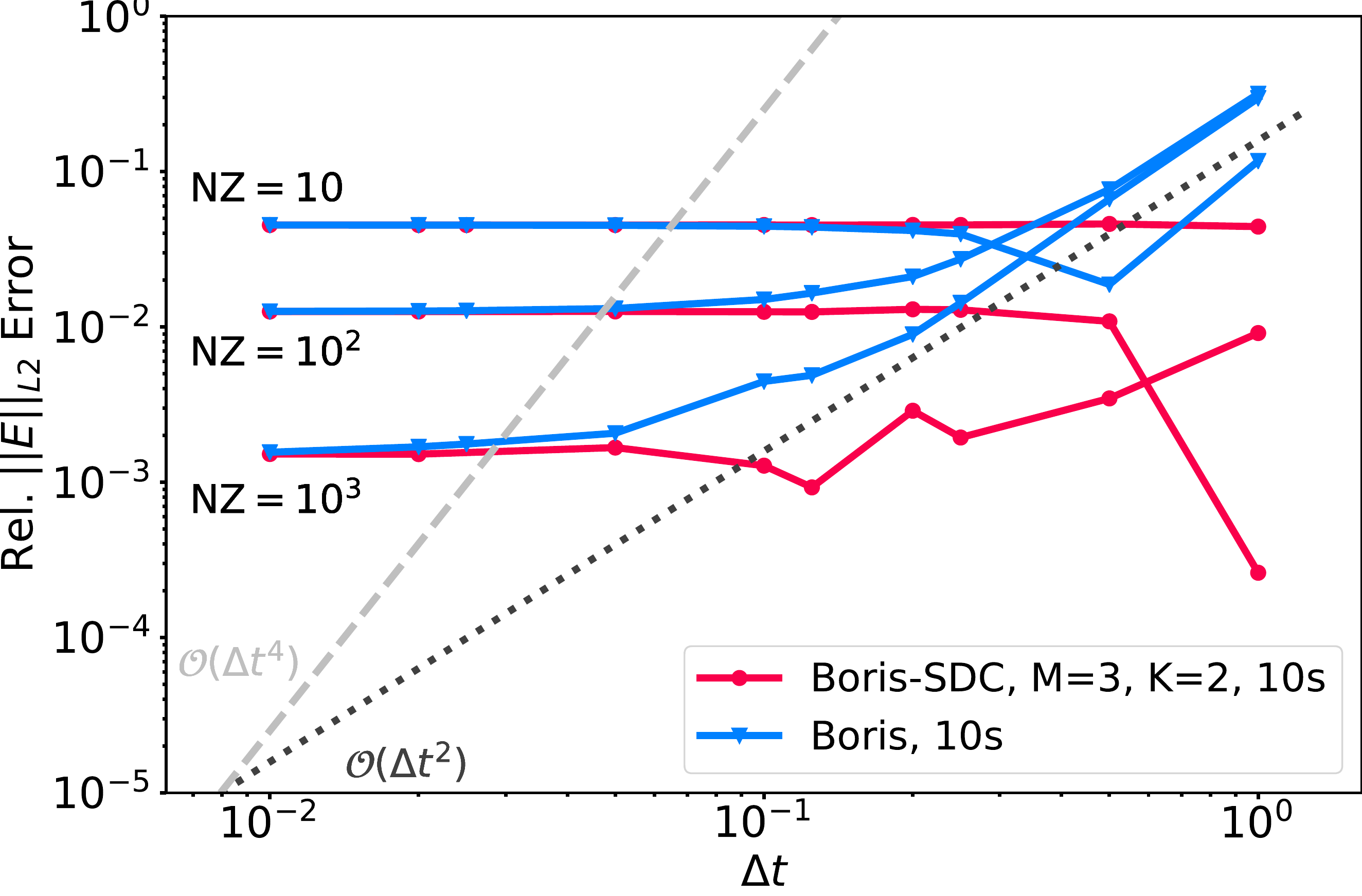}
    \end{minipage}
    \begin{minipage}{0.475\textwidth}
        \centering
        \includegraphics[width=1\linewidth]{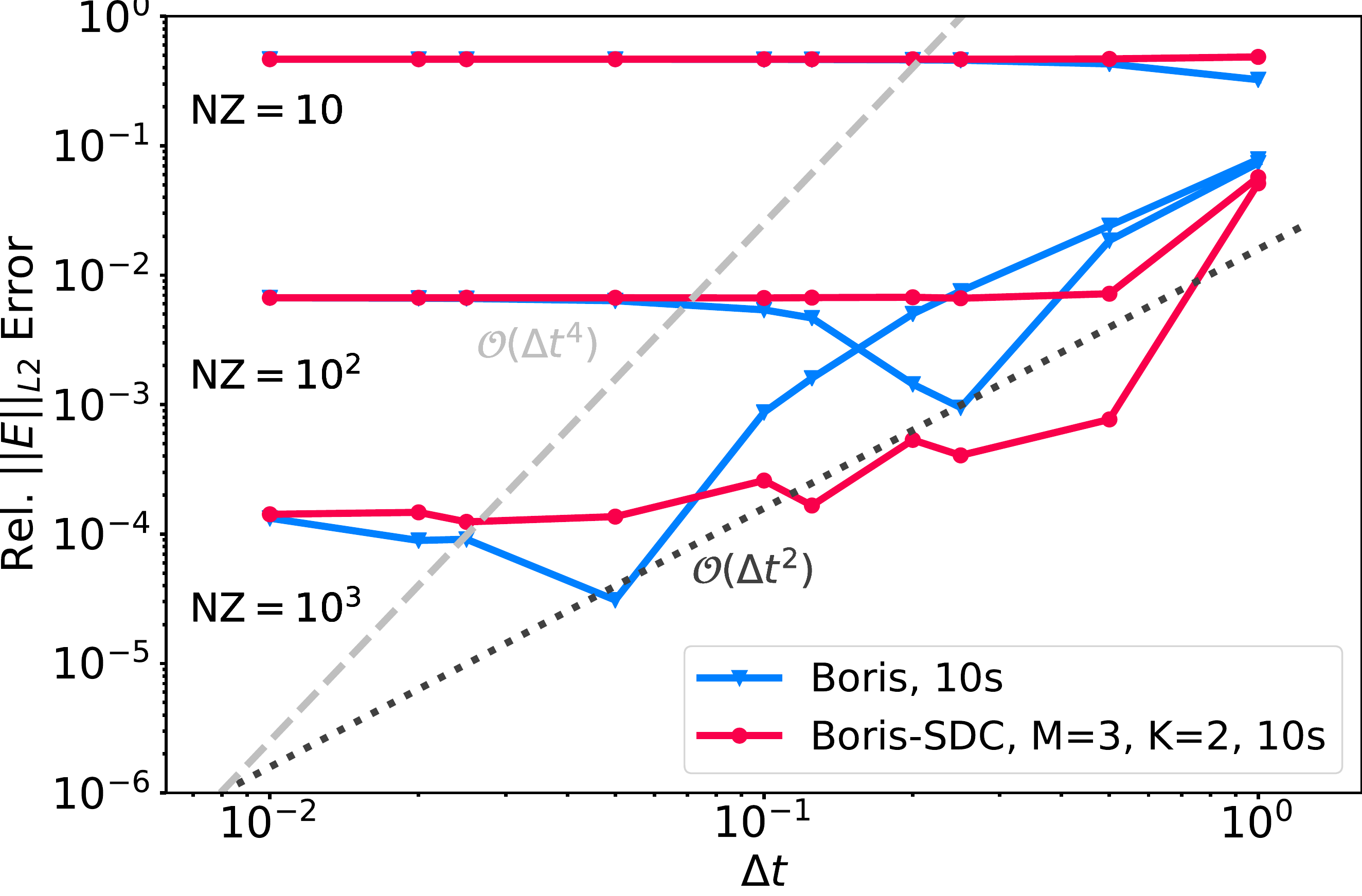}
     \end{minipage}
    \caption{Error at \SI{10}{\second} versus time step size for weak (left) and strong (right) Landau damping.}
    \label{fig:landau_dt}
\end{figure}

Figure~\ref{fig:landau_dt} shows the error in electric field norm compared to the reference simulation at $t=10$ against time step size for weak (left) and strong (right) Landau damping.
In most cases, the spatial error dominates and there is little effect from varying $\Delta t$.
However, for $N_z = 10^3$, Boris-SDC is more accurate for the same $\Delta t$ and reaches the saturation error set by the spatial resolution earlier than the Boris method.
For example, Boris-SDC reaches an error of $10^{-3}$ for weak Landau damping with a step size of around $\Delta t = 10^{-1}$, compared to the Boris integrator which requires $\Delta t = 10^{-2}$.
This effect is more pronounced for the weak Landau damping, most likely because the smoother dynamics lead to a smaller spatial discretisation error.

\begin{figure}[t]
    \begin{minipage}[c]{0.475\textwidth}%
        \centering
        \includegraphics[width=0.95\linewidth]{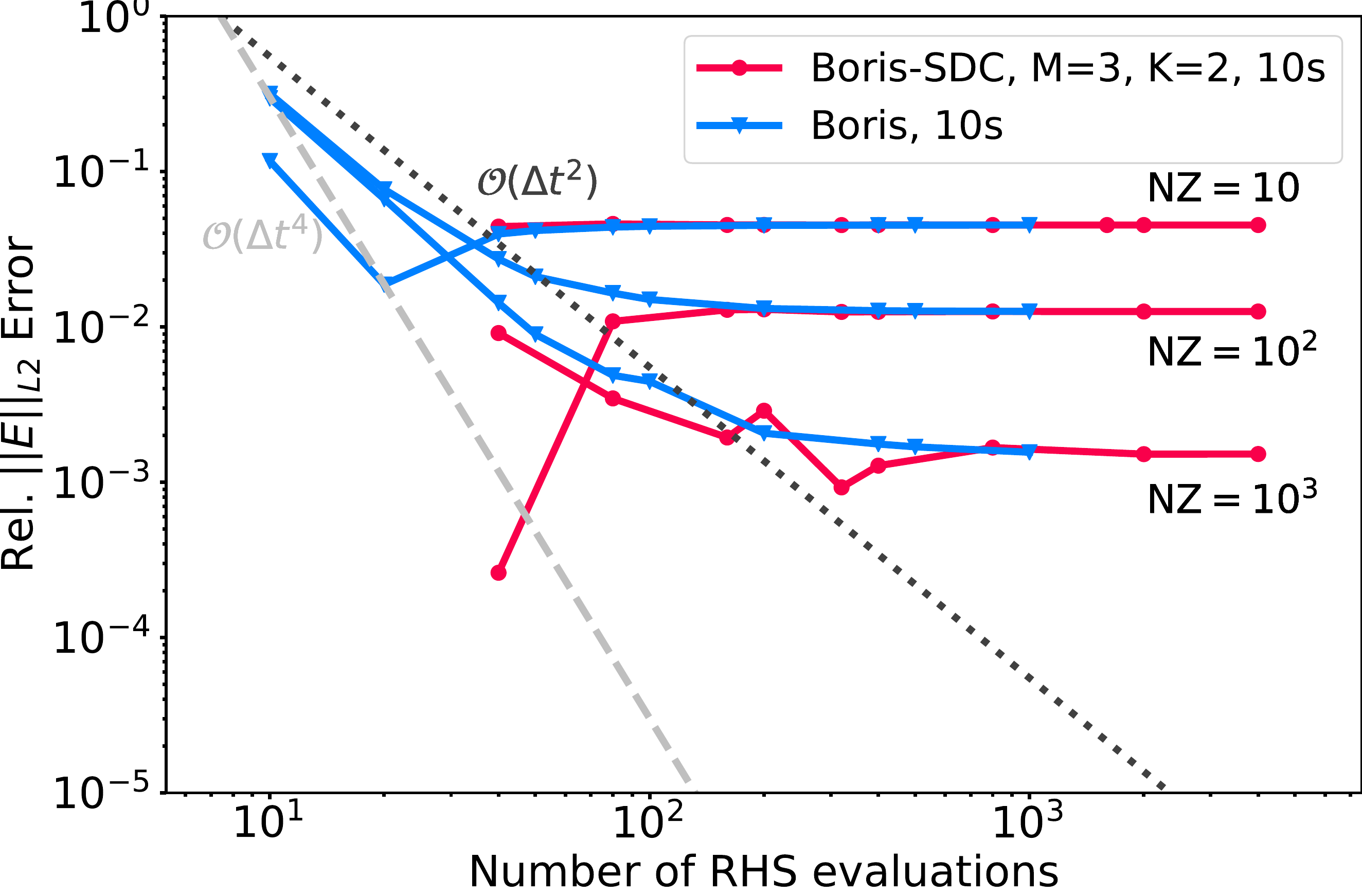}
    \end{minipage}
    \begin{minipage}[c]{0.475\textwidth}%
        \centering
        \includegraphics[width=0.95\linewidth]{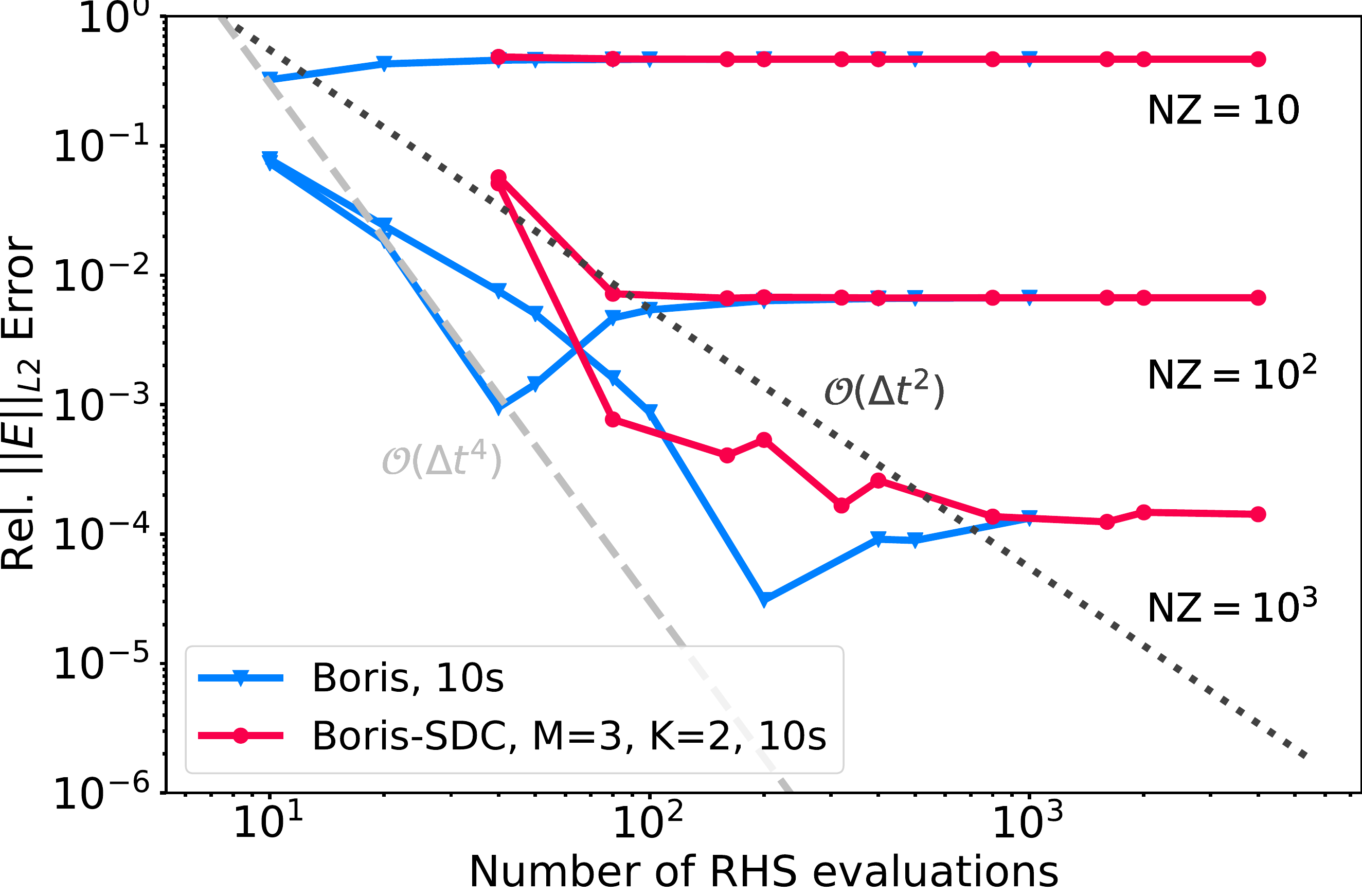}
    \end{minipage}
    \caption{Error versus number of right hand side evaluations for weak (left) and strong (right) Landau damping.}
    \label{fig:landau_rhs}
\end{figure}

Figure~\ref{fig:landau_rhs} shows the error at $t=10$ against computational effort, measured by the number of right hand side evaluations. 
For weak Landau damping, both methods perform similarly.
Although Boris-SDC allows one to achieve a given accuracy with a larger time step size, the reduced computational effort from computing fewer time steps is counterbalanced by the increased workload per step.
Only minimal gains are achieved for weak Landau damping for errors between $10^{-2}$ and $10^{-3}$ where Boris-SDC is marginally more efficient.
For strong Landau damping, we do not see efficiency gains from Boris-SDC, despite its better accuracy.

\subsection{Relativistic Penning Trap}
\begin{figure}[t]
	\centering
	\begin{minipage}[c]{0.45\textwidth}%
	\centering
	\includegraphics[width=1\linewidth]{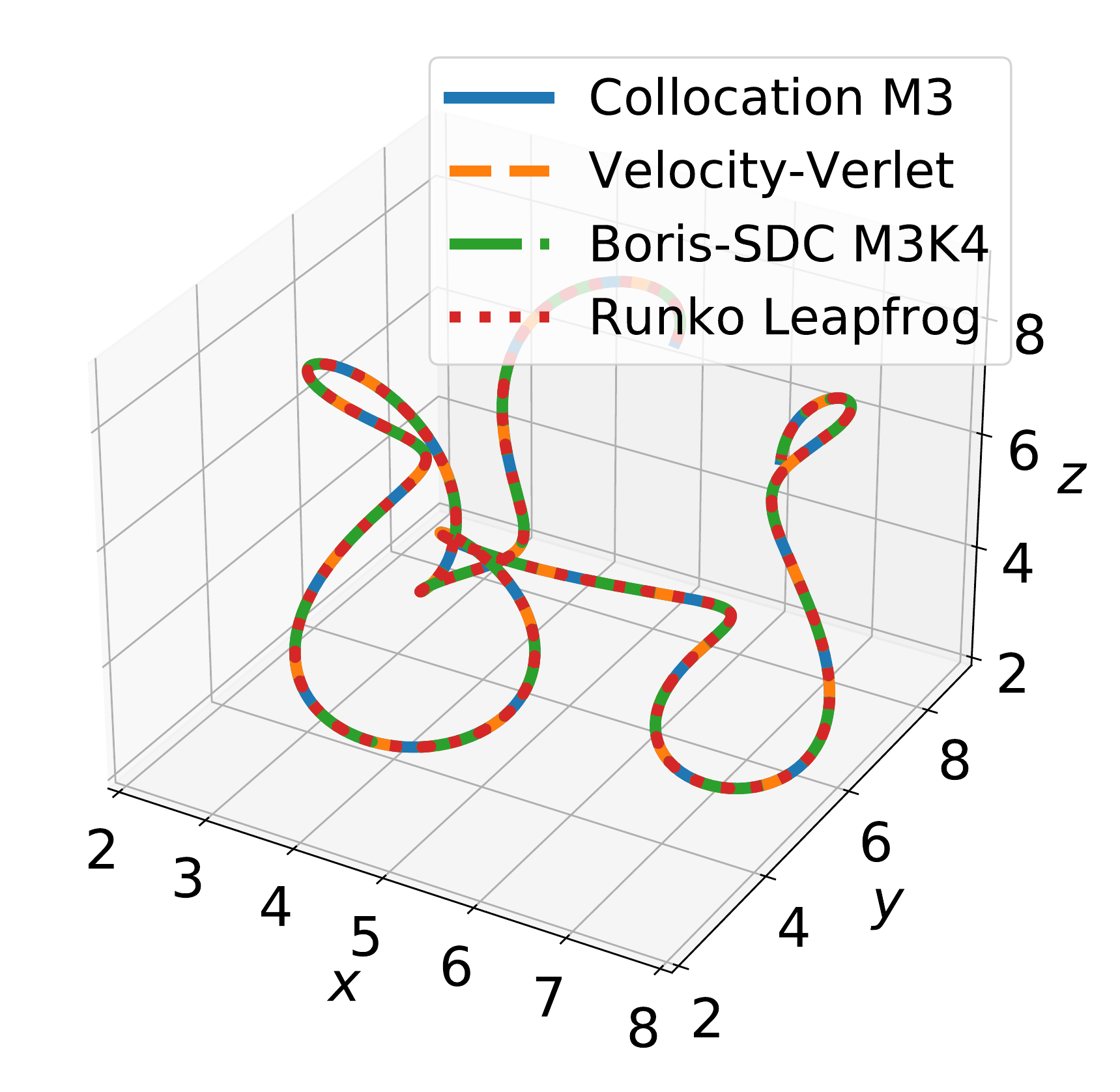}
	\caption{Runko and Boris-SDC M3K4 Penning trajectory for $t=[0,45]$.}
	\label{fig:runko_penning_iso}
	\end{minipage}\hspace*{\fill}
	\begin{minipage}[c]{0.45\textwidth}%
	\centering
	\includegraphics[width=1\linewidth]{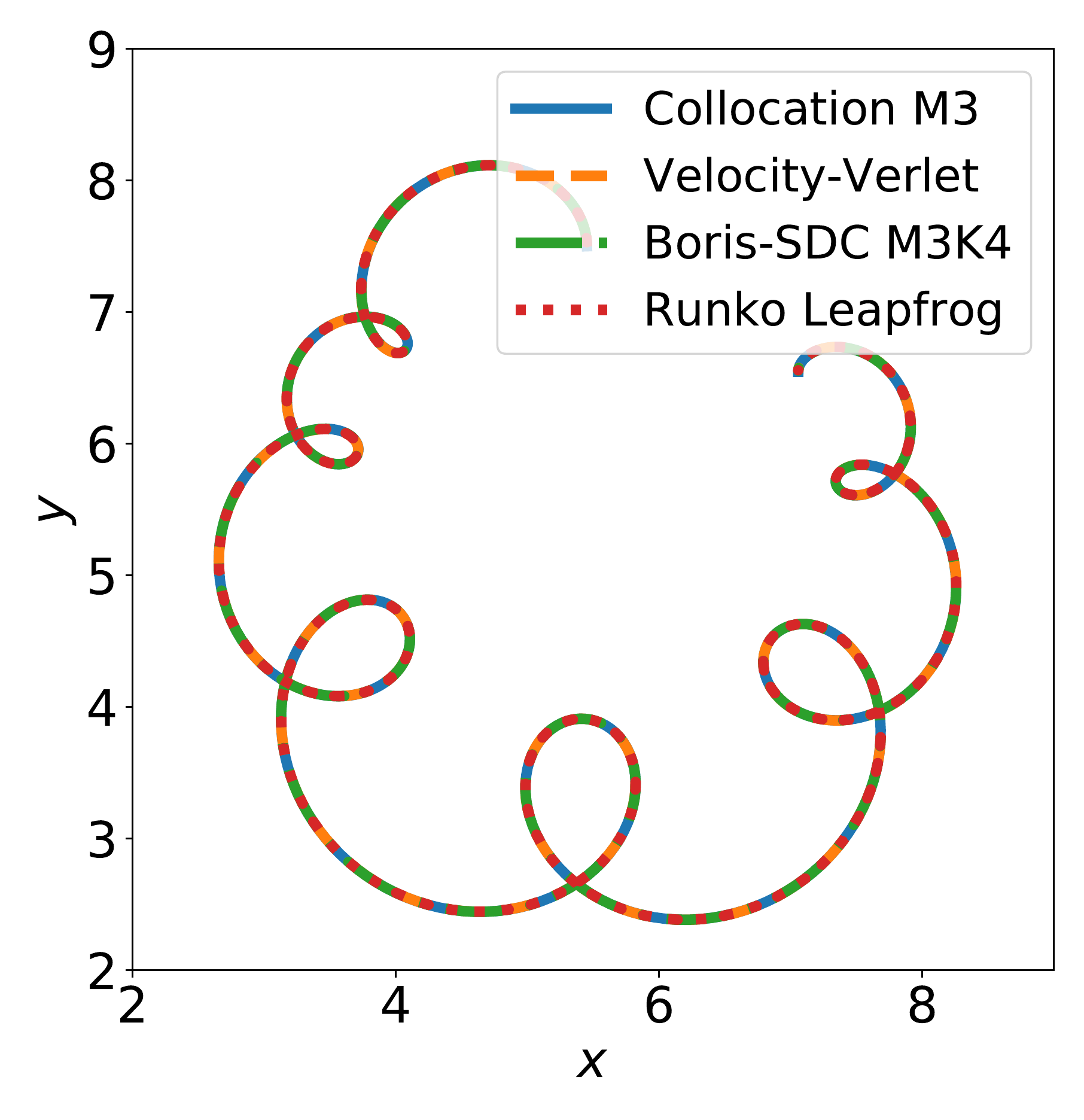}
	\caption{Runko and Boris-SDC M3K4 Penning xy-plane trajectory for $t=[0,45]$.}
	\label{fig:runko_penning_xy}
	\end{minipage}
\end{figure}

\begin{figure}[t]
    \centering
    \includegraphics[width=.49\linewidth]{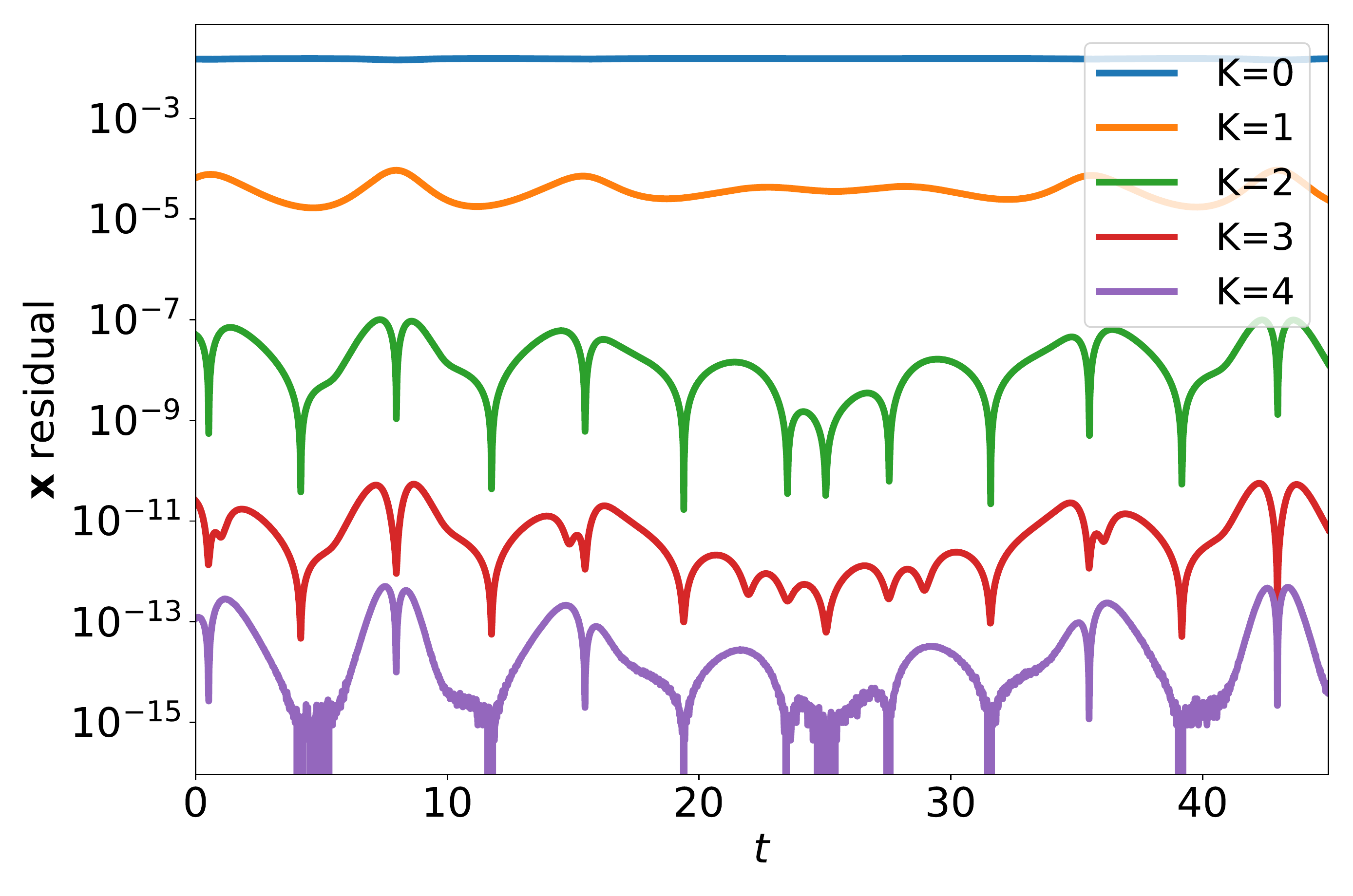}
    \includegraphics[width=.49\linewidth]{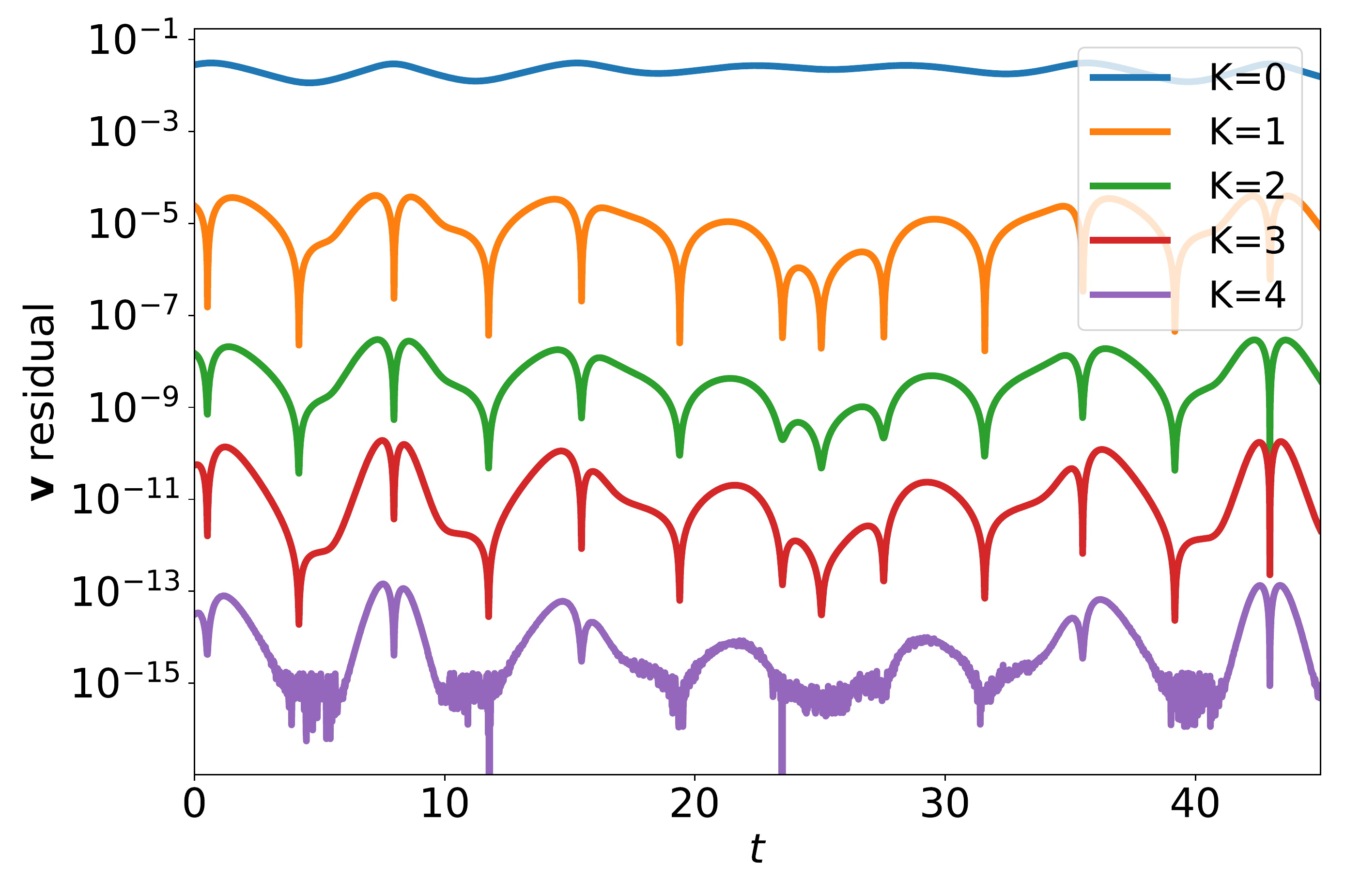}
    \caption{Boris-SDC residual in position (left) and velocity (right) for the relativistic Penning trap.}
    \label{fig:runko_penning_res}
\end{figure}

We compare performance of Boris and Boris-SDC in terms of work-precision for the Penning trap, similar to the non-relativistic test cased used by Winkel et al.~\cite{winkel2015highOrderBoris}.
Here, however, we translate the parameters into the units used by Runko.
The simulation parameters are summarised in Table~\ref{tab:runko_penning_params}.
The fields at a grid node with index $(i,j,k)$ are
\begin{equation}
E_{i,j,k}=|E|\begin{bmatrix}1 & 0 & 0\\
0 & 1 & 0\\
0 & 0 & -2
\end{bmatrix}x_{i,j,k},\label{eq:runko_penning_E}
\end{equation}
and
\begin{equation}
B_{i,j,k}=|B|,\label{eq:runko_penning_B}
\end{equation}
where $x_{i,j,k}$ is the position of node.
The result is a homogeneous magnetic field pointing up along the $z$-axis and an electric field pushing towards and along the $xy$-plane centred on $z=0$. Charges caught in the fields are pushed towards this plane and away from the centre by the electric field, with the magnetic field curving the trajectories back to keep the charges trapped.
We use the linear interpolation to compute field values at particle positions and the Boris integrator as a reference particle pusher. 
The resulting trajectory and trajectory projection on the $x$-$y$-plane can be seen in Figures~\ref{fig:runko_penning_iso} and~\ref{fig:runko_penning_xy}.
\begin{table}[t]
	\centering{}
	\caption{Validation study parameters}
	\label{tab:runko_penning_params} %
	\begin{tabular}{cc}
	\hline
	$\hat{c}$ & 0.45\tabularnewline
	$t_{\mathrm{end}}$ & 45\tabularnewline
	$\hat{\mathbf{x}}(0)$ & $(7.5,5,7.5)^{\mathrm{T}}$\tabularnewline
	$\hat{\mathbf{v}}(0)$ & $(0.7\hat{c},0,0.7\hat{c})^{\mathrm{T}}$\tabularnewline
	$N_x, N_y, N_z$ & $10$\tabularnewline
	$|E|$ & 0.1\tabularnewline
	$|B|$ & 1\tabularnewline
	\hline
	\end{tabular}
\end{table}
Figures~\ref{fig:runko_penning_res} shows the SDC residual in position (left) and velocity (right) for Boris-SDC for $K=0$ up to $K=4$ iterations.
Because of the choice of the Lorentz parameter described in Subsection~\ref{sec:3-rel-boris}, the Lorentz factor $\gamma$ converges to the correct values as $K$ increases and the residual goes down to more or less machine precision as the iteration converges to the collocation solution.

\begin{figure}[t]
	\centering
	\includegraphics[width=.49\linewidth]{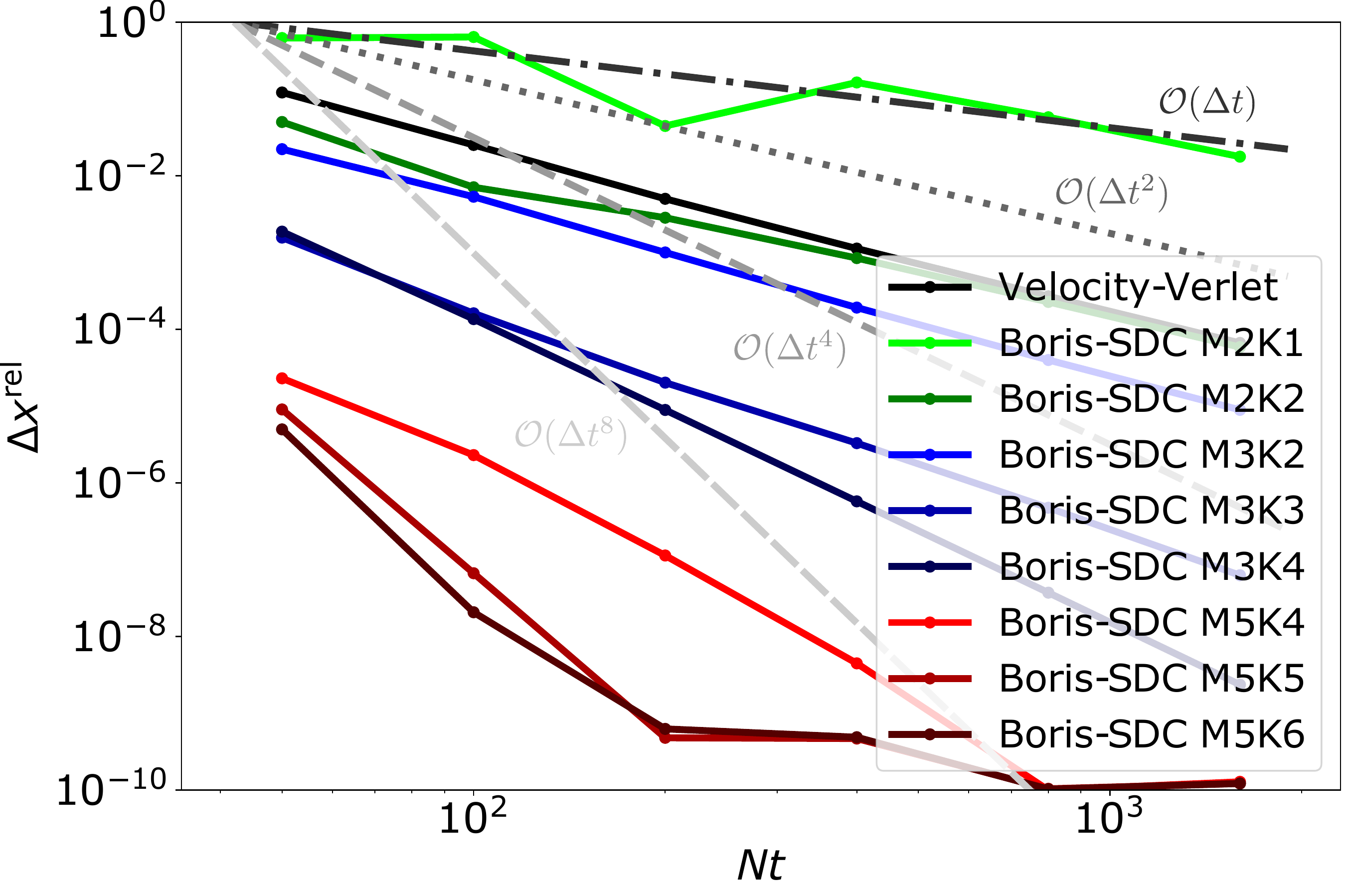}~~~~
	\includegraphics[width=.49\linewidth]{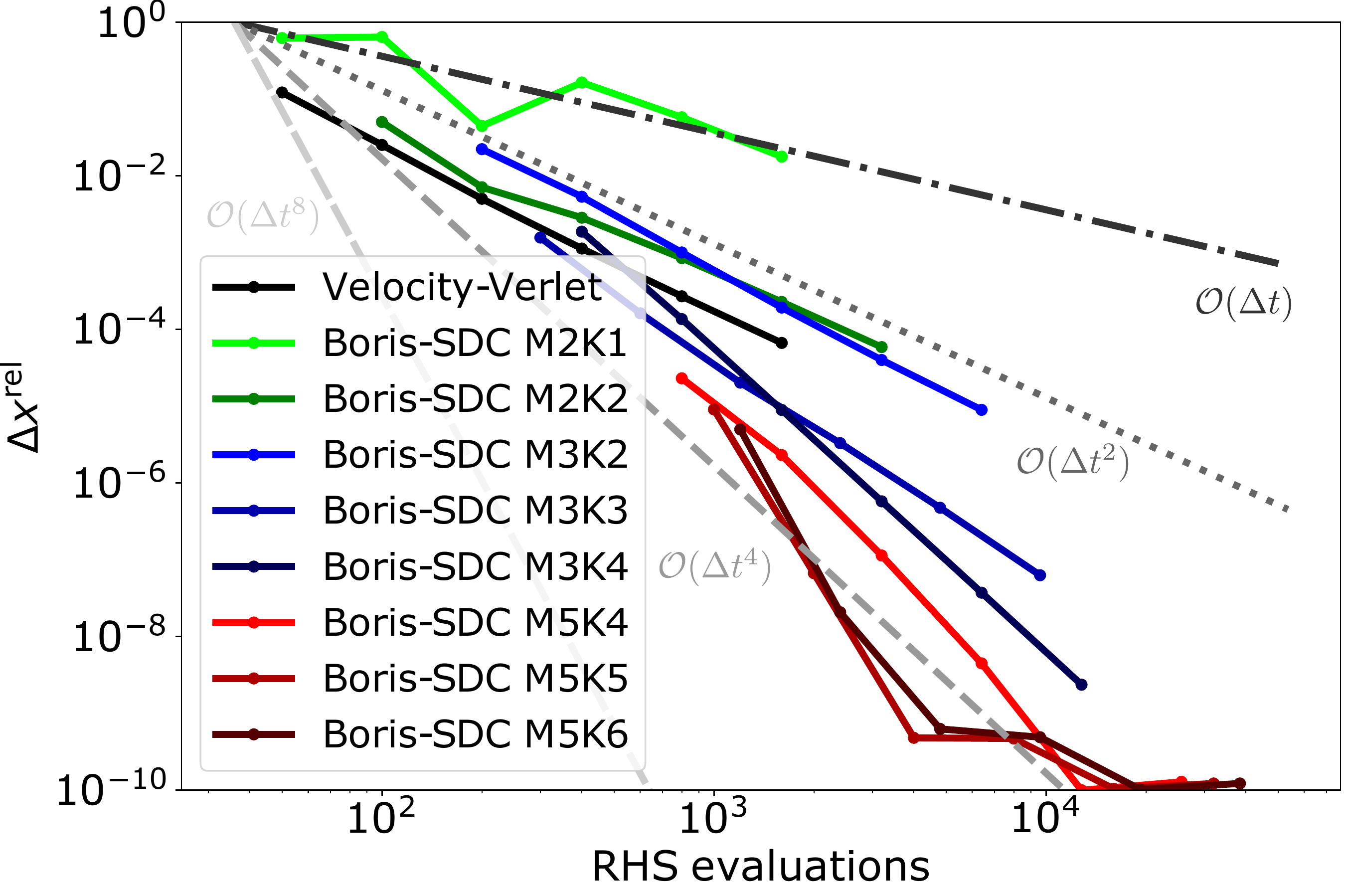}
	\caption{Error versus number of time steps (left) and error versus computational cost (right) for the relativistic Penning trap.}
	\label{fig:runko_penning_wp}
\end{figure}

To compare Boris and Boris-SDC in terms of work-precision, we compute a reference solution using the same $10 \times 10 \times 10$ grid, $NT = 3200$ and $M=5$ quadrature nodes.
Figure~\ref{fig:runko_penning_wp} shows error against time step (left) and error against computational cost (right), measured again by the number of required right hand side evaluations.
As a guide to the eye, gray lines with slopes of minus one, two, four and eight are also shown.
As expected, Boris/velocity-Verlet converges with order two.
For Boris-SDC, the order increases with $K$, although the precise impact of an additional iteration is not clear.
While Winkel et al.~\cite{winkel2015highOrderBoris} observed numerically for the non-relativistic case that a sweep increases order by two (until reaching the order of the collocation method), the picture is less clear in the relativistic case.
For both $M=3$ and $M=5$ nodes, Boris-SDC reproduces the order four or eight, respectively, of the underlying collocation method for sufficiently many iterations.
Furthermore, every iteration leads to a significant gain in accuracy, even when it fails to fully increase the order by two.

\subsection{Numerical drift in the force-free case}
If the acceleration from the magnetic and electric field in~\eqref{eq:lorentz_rel} cancel, there is no net force acting on the particle and it should continue to travel into its original direction without changing velocity.
However, owing to round-off error, acceleration from the fields will not cancel out exactly on the discrete level and the particle will undergo numerical drift.
For the test case studied by Ripperda et al.~\cite{ripperda2018pushers}, we will compare numerical drift for the Boris integrator, the Vay integrator~\cite{vay2008simulation} and Boris-SDC with $M=3$ and $M=5$ nodes and varying iteration counts.
The test cases uses a particle with an initial velocity that is parallel to the $y$-axis and a magnetic field with field lines oriented along the $z$-axis.
Then, the electric field is set to
\begin{equation}
	\mathbf{E} = - \mathbf{v} \times \mathbf{B} = \begin{pmatrix} v_y B_z  \\ 0 \\ 0 \end{pmatrix} = \begin{pmatrix} E_x \\ 0 \\ 0 \end{pmatrix}
\end{equation}
so that the net acceleration is zero.
Parameters for the problem are summarised in Table~\ref{tab:force_free_par}.
\begin{table}[]
    \centering
    \caption{Force-free test parameters}
    \label{tab:force_free_par} %
    \begin{tabular}{ccccc}
    \toprule 
    Parameter & Value \tabularnewline
    \midrule
    $c$ & $1$ \tabularnewline
    $q/m$ & $1$ \tabularnewline
    $\gamma$ & $10^6$ \tabularnewline
    $B_z$ & $1$ \tabularnewline
    $v_y$ & $\sqrt{1-\frac{1}{\gamma^2}} \cdot c$\tabularnewline
    $E_x$ & $-v_x B_z$\tabularnewline
    $T_{end}$ & $10^5$\tabularnewline
    \bottomrule
    \end{tabular}
\end{table}

\begin{figure}[]
\centering
\begin{minipage}[c]{0.49\textwidth}%
\includegraphics[width=1\linewidth]{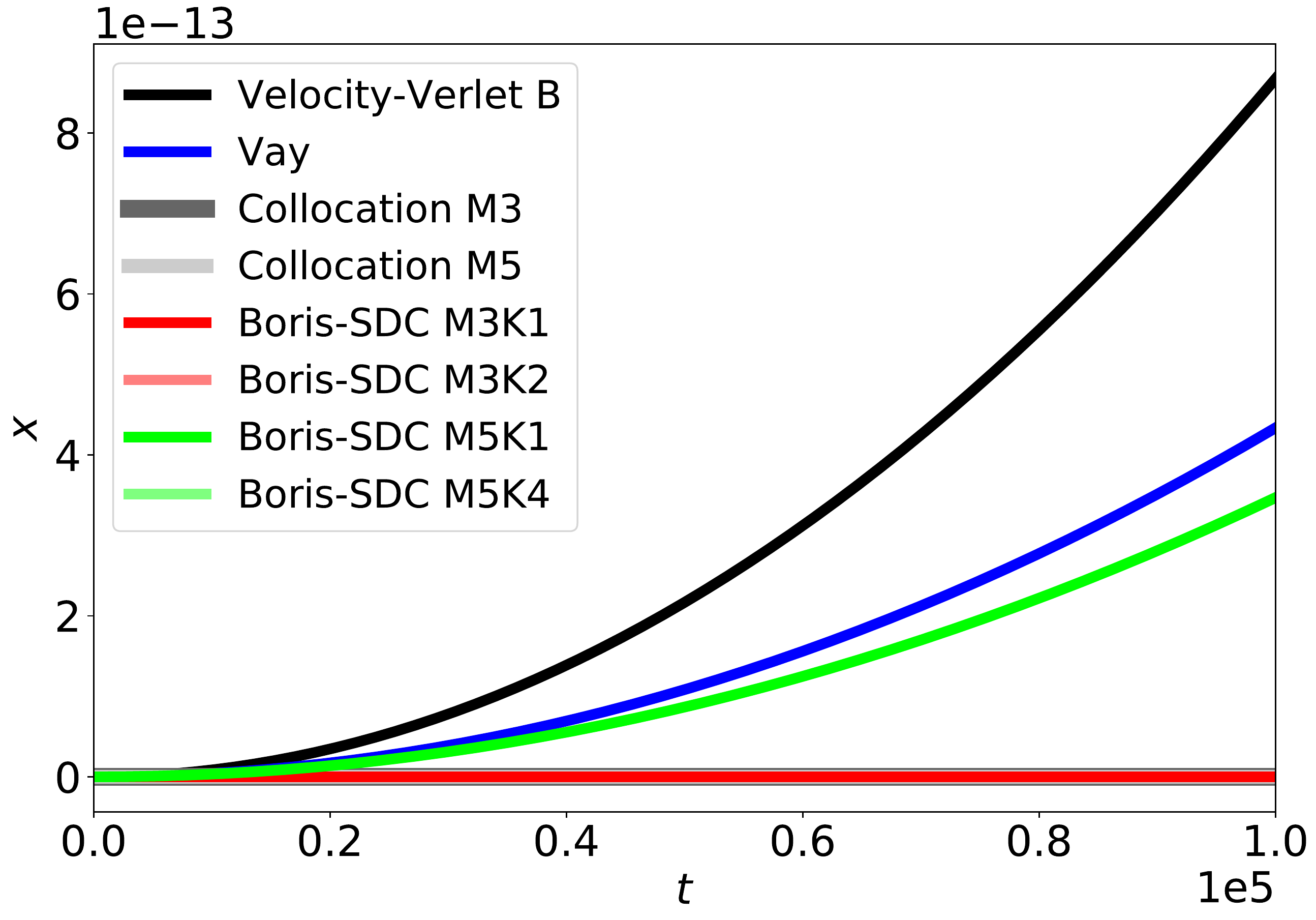}
\end{minipage}\hspace*{\fill}%
\begin{minipage}[c]{0.51\textwidth}%
\includegraphics[width=1\linewidth]{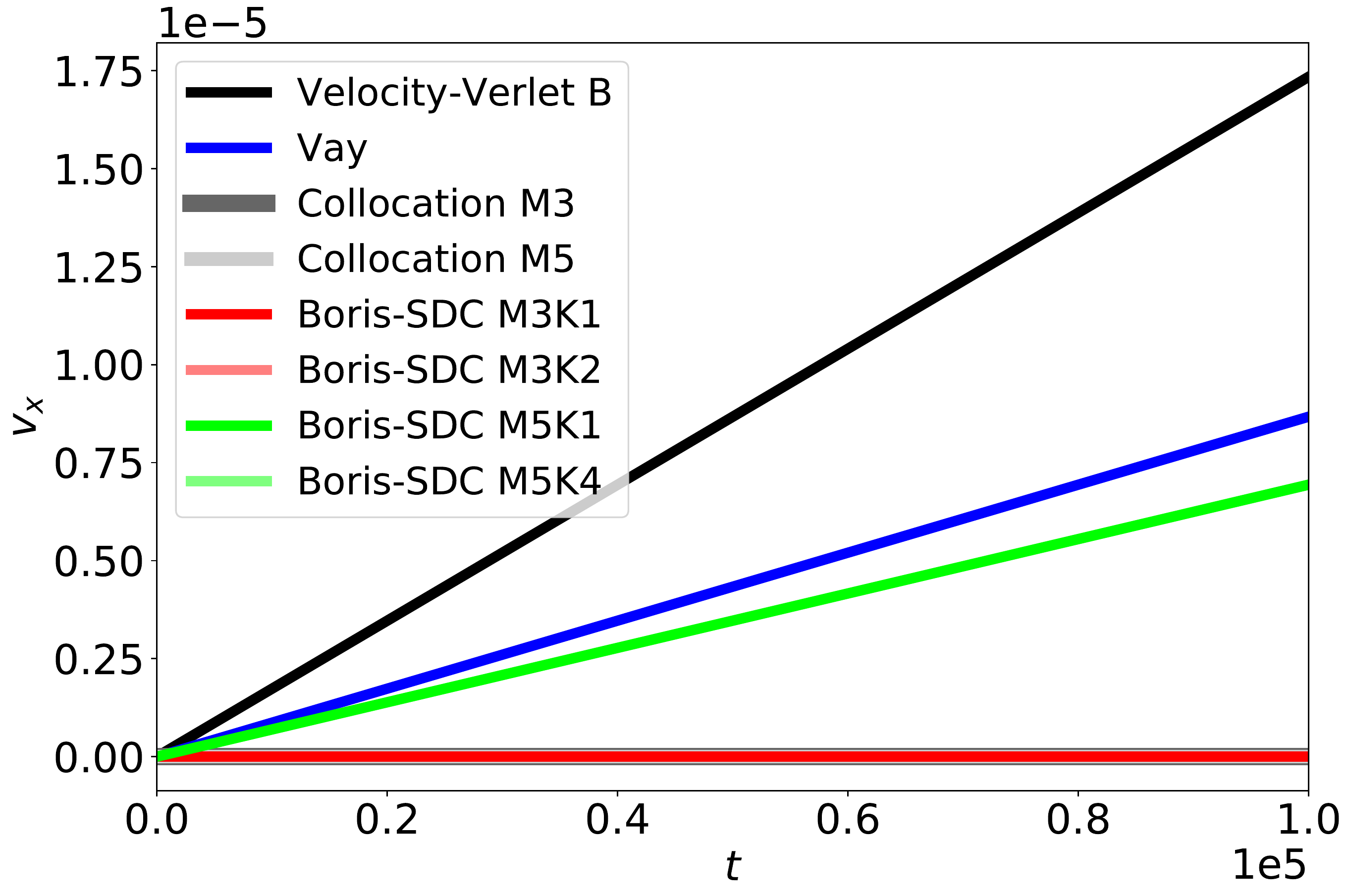}
\end{minipage}
\caption{Error in position (left) and velocity (right) for the force-free case. Note that the lines for Collocation M3 and M5 and Boris-SDC M3K2 and M5K4 all overlap.}
\label{fig:force_free_drift}
\end{figure}

Figure~\ref{fig:force_free_drift} shows the resulting drift in position (left) and velocity (right) over time for Boris, Vay and Boris-SDC with $M=3$ nodes and $K=1, 2$ iterations and $M=5$ nodes and $K=1,4$ iterations.
For both $M=3$ and $M=5$ nodes, the collocation solution exhibits no numerical drift.
Boris-SDC with $M=3$ nodes is drift-free already after a single iteration.
When $M=5$ nodes are used, Boris-SDC with $K=1$ iteration suffers from drift, most likely due to floating point arithmetic imprecision.
However, its drift is already less than that of the Vay method and much less than that of classical Boris. Across all studies of the force-free drift, Boris-SDC exhibits either no drift or drift comparable to that of Vay, with none of the exponentially scaling drift at coarse time-steps seen for classical Boris.

%% file: conclusions.tex
\section{Conclusions and Outlook}\label{sec:5-conclusions}
The paper explores the effect of using Boris spectral deferred corrections (Boris-SDC), a high order generalisation of the Boris algorithm, as particle pusher in a particle-in-cell code.
Accuracy and performance in terms of work-precision is investigated for simulations of a two-stream instability and Landau damping.
A modification for Boris-SDC is proposed to make it applicable to the relativistic Lorentz equations.
Performance of relativistic Boris-SDC is then studied for a relativistic Penning trap and numerical drift is assessed for a particle in a scenario where accelerations from the electric and magnetic field should cancel out.

Throughout, we found that Boris-SDC is more accurate than Boris at the same time step size, yielding lower error for any given time-step size, up to the limits imposed by the accuracy of the spatial discretisation.
Alternatively, a given accuracy could be reached with a much larger time step compared to Boris.
However, for the non-relativistic cases, these gains were offset by the much higher computational cost per time step; as a result, Boris and Boris-SDC performed very similarly in terms of overall computational cost.
A key problem was that mostly errors from the spatial discretisation were dominating, limiting the gains from the higher order of Boris-SDC.
Results were somewhat more positive for the relativistic case, partly because the imposed electric and magnetic fields are relatively simple and can be interpolated exactly.
We confirmed that relativistic Boris-SDC delivers higher order of convergence, up to the order of the underlying quadrature method.
This also translated into computational gains compared to relativistic Boris for errors below $10^{-3}$.
Finally, tests for the force-free case showed that Boris-SDC produces less numerical drift than both the Boris method and Vay integrator.
Some configurations of Boris-SDC were even completely drift free, but we do not have, at the moment, a theoretical understanding why this is.

Although the results so far do not show clear efficiency gains from Boris-SDC, we expect that PIC/Boris-SDC could outperform PIC/Boris when combined with high order spatial approximations as well as techniques to further improve performance of SDC~\cite{EmmettMinion2012,HuangEtAl2006,Weiser2014} and in situations where high accuracy is required.
Potential practical applications include, for example, studies of growth rate and saturation of especially weak plasma instabilities in laboratory and astrophysical context.
Topics of interest for further exploration  include more complex relativistic test cases as well as the development of a fully electromagnetic Boris-SDC/PIC scheme. 
Another promising avenue of research would be to understand the surprising drift-free nature of some configurations of Boris-SDC: it might be possible to take advantage of this feature by applying it for example to magnetized collisionless shock simulations where high numerical precision is needed~\cite{SironiEtAl2021}.